\documentclass[english,x11names]{amsart}

\PassOptionsToPackage{hyphens}{url}
\usepackage{amsxtra}
\usepackage{mathtools}
\usepackage{enumitem}
\usepackage{amsfonts}
\usepackage[toc,page]{appendix}
\usepackage{amsthm}
\usepackage[most]{tcolorbox}
\usepackage{graphicx}
\usepackage{amssymb}
\usepackage{extarrows}
\usepackage[protrusion=true, tracking=true, expansion=true, final]{microtype}
\usepackage{tikz-cd}
\usepackage{tikz}
\usetikzlibrary{calc}
\usepackage{textcomp}
\usepackage{fontawesome}
\usepackage{hyphenat}
\usepackage{mathabx}
\usepackage[bbgreekl]{mathbbol}
\usepackage{relsize}
\usepackage[T1]{fontenc}
\usepackage{charter}
\usepackage{seqsplit}
\usepackage{breqn}
\usepackage{upgreek}
\usepackage{mathrsfs}
\usepackage[colorlinks=true, citecolor=DodgerBlue2, linkcolor=HotPink2, urlcolor=Aquamarine3]{hyperref}

\theoremstyle{plain}
\newtheorem{thm}{Theorem}[section]
\newtheorem{cor}[thm]{Corollary}
\newtheorem{pro}[thm]{Proposition}
\newtheorem{lem}[thm]{Lemma}
\theoremstyle{definition}
\newtheorem{dfn}[thm]{Definition}
\newtheorem{example}[thm]{Example}
\newtheorem{rmk}[thm]{Remark}
\newtheorem{rem}[thm]{Remark}

\newtheorem*{conjecture*}{Conjecture}
\newtheorem{que}[thm]{Question}
\newtheorem{construction}[thm]{Construction}
\newtheorem*{assumption*}{Assumption}

\newtheorem{setup}[thm]{Setup}

\newtheorem{cond}[thm]{Condition}

\DeclareSymbolFontAlphabet{\mathbb}{AMSb}
\DeclareSymbolFontAlphabet{\mathbbl}{bbold}
\newcommand{\prism}{{\mathlarger{\mathbbl{\Delta}}}}

\def\into{\hookrightarrow}

\def\Hom{\mathop{\textup{Hom}}\nolimits}

\def\Gal{\mathop{\textup{Gal}}\nolimits}
\def\End{\mathop{\textup{End}}\nolimits}

\def\Perf{\mathop{\bf Perf}\nolimits}

\def\perf{\mathop{\textup{perf}}\nolimits}

\def\Spec{\mathop{\bf Spec}\nolimits}
\def\Spa{\mathop{\bf Spa}\nolimits}
\def\Spf{\mathop{\bf Spf}\nolimits}
\def\Spd{\mathop{\bf Spd}\nolimits}
\def\Sht{\mathop{\bf Sht}\nolimits}
\def\Sh{\mathop{\bf Sh}\nolimits}
\def\Bun{\mathop{\bf Bun}\nolimits}
\def\Div{\mathop{\bf Div}\nolimits}

\def\kernel{\mathop{\textup{Ker}}\nolimits}

\def\GL{\mathop{\textup{GL}}\nolimits}
\def\Gr{\mathop{\textup{Gr}}\nolimits}

\def\red{{\textup{red}}}
\def\ad{{\textup{ad}}}
\def\der{{\textup{der}}}

\def\Fib{{\textup{Fib}}}
\def\der{{\textup{der}}}

\def\id{{\textup{id}}}
\def\Flex{{\textup{Flex}}}
\def\com{{\textup{com}}}

\def\Res{{\textup{Res}}}

\newcommand{\BA}{{\mathbb{A}}}

\newcommand{\BC}{{\mathbb{C}}}

\newcommand{\BF}{{\mathbb{F}}}
\newcommand{\BG}{{\mathbb{G}}}

\newcommand{\BN}{{\mathbb{N}}}
\newcommand{\BO}{{\mathbb{O}}}

\newcommand{\BQ}{{\mathbb{Q}}}
\newcommand{\BR}{{\mathbb{R}}}

\newcommand{\BZ}{{\mathbb{Z}}}

\newcommand{\FX}{{\mathfrak{X}}}

\newcommand{\CB}{{\mathcal{B}}}

\newcommand{\CD}{{\mathcal{D}}}
\newcommand{\CE}{{\mathcal{E}}}
\newcommand{\CF}{{\mathcal{F}}}
\newcommand{\CG}{{\mathcal{G}}}

\newcommand{\CK}{{\mathcal{K}}}

\newcommand{\CM}{{\mathcal{M}}}
\newcommand{\CN}{{\mathcal{N}}}
\newcommand{\CO}{{\mathcal{O}}}

\newcommand{\CS}{{\mathcal{S}}}

\newcommand{\CY}{{\mathcal{Y}}}

\newcommand{\SD}{{\mathscr D}}

\newcommand{\SM}{{\mathscr M}}

\newcommand{\SP}{{\mathscr P}}
\newcommand{\SQ}{{\mathscr Q}}

\newcommand{\xqed}[1]{%
  \leavevmode\unskip\penalty9999 \hbox{}\nobreak\hfill
  \quad\hbox{\ensuremath{#1}}}
  
\title{Embeddings of Certain Exceptional Shimura Varieties into Siegel Modular Varieties}

\author{S. Mohammad Hadi Hedayatzadeh}
\address{Institute for Research in Fundamental Sciences (IPM),
Niavaran Sq., Tehran, Iran}
\email{hadi@ipm.ir}

\author{Ali Partofard}
\address{Institute for Research in Fundamental Sciences (IPM),
Niavaran Sq., Tehran, Iran}
\email{alipartofard@ipm.ir}

\begin{document}

\begin{abstract}
    We define a class of local Shimura varieties that contains some local Shimura varieties for exceptional groups, and for this class, we construct a functor from $\left(G, \mu\right)$-displays to $p$-divisible groups. As an application, we prove that for this class, the local Shimura variety is representable and perfectoid at the infinite level. Considering the global counterpart of this class, we embed certain exceptional Shimura varieties into Siegel modular varieties. In particular, we prove that they are perfectoid at the infinite level.
\end{abstract}
\maketitle

\begin{center}
    \noindent\rule{0.8\linewidth}{0.4pt}\vspace{5pt}

    \small{2020 Mathematics Subject Classification: 14F30 (Primary); 14L05 (Secondary)}

    \noindent\rule{0.8\linewidth}{0.4pt}
\end{center}
\vspace{0.5cm}

\tableofcontents


\setcounter{section}{0}
\phantomsection
\section*{Introduction}


Let $(G,\mu)$ be a \emph{Shimura datum}: $G$ is a reductive group over $\mathbb{Q}$, and $\mu$ is a cocharacter; they satisfy certain conditions. One can attach to this data a tower of smooth varieties defined over $\mathbb{C}$, called \emph{Shimura variety}, and one can show that this tower has a canonical model over a number field. The theory of Shimura varieties plays a central role in the Langlands program, providing the geometric framework for relating automorphic representations to Galois representations. The simplest examples of Shimura varieties are the Siegel modular varieties, the moduli space of polarized Abelian varieties with level structures.\\

Using this moduli interpretation, many nice results have been proved for Siegel modular varieties. For example, we know that at parahoric levels, there are canonical integral models for the Siegel modular variety that satisfy a Néron-like extension property. There is a group-theoretic description of the points of the Siegel modular variety over finite fields; this helps us compute the Hasse-Weil $\zeta$-function associated with this variety. Scholze \cite{scholze2015torsion} proved that at the level $p^\infty$ the Siegel modular variety is representable by a perfectoid space.\\

While Siegel modular varieties serve as the archetype for the theory of Shimura varieties, many fundamental properties established for them, which are central to the Langlands program, are expected to hold for all Shimura varieties. In fact, according to Langlands' \emph{Märchen}, all Shimura varieties should behave as if they are moduli spaces of motives with additional structures \cite{langlands1979marchen}. A Shimura datum is said to be of \emph{Abelian type} if, after a central isogeny, it can be embedded into the Siegel datum. Using this embedding, sometimes called ``Deligne’s trick'', one can generalize all the above results to the case of Shimura varieties of Abelian type. Through the work of Kisin and others, results such as the existence of canonical integral models and the Langlands-Rapoport conjecture have been extended to Shimura varieties of Abelian type \cite{kisin2010integral, kisin2017mod}. Xu Shen \cite{shen2017perfectoid} used Deligne's formalism to prove that every Shimura variety of Abelian type with infinite level at $p$ is perfectoid.\\

The local counterpart to Shimura varieties, the so-called \emph{local Shimura varieties}, conjectured by Rapoport and Viehmann \cite{RapoportViehmann2014}, was defined by Scholze, who also proved their representability by rigid spaces \cite{scholze2020berkeley}. Scholze also constructed their integral model as a v-stack. Similar to the global setting, the prototypes for the local Shimura varieties are Rapoport-Zink spaces: the moduli spaces of $p$-divisible groups. Again, one can prove many nice results about Rapoport-Zink spaces and, more generally, local Shimura varieties of local Abelian type. For example, one can prove the representability of the integral model \cite{pappas2022p} and the fact that they are representable by perfectoid spaces at the $p^\infty$-level.\\

For Shimura varieties of non-(pre-)Abelian type, such as those associated with the exceptional group $E_7$, the above strategy fails, as the associated Hodge structures no longer correspond to Abelian varieties, and in the local case, the associated crystals do not correspond to $p$-divisible groups.\\

Progress has been made in extending these local and global results to general Shimura data; yet, significant questions still remain. For instance, T. He \cite{he2024perfectoidness} proved the perfectoidness of completed stalks for general Shimura varieties. However, the perfectoidness of Shimura varieties has not been established in general. In the local setting, while Gleason \cite{gleason2020specialization} proved the representability of the local model for any local Shimura variety with parahoric level structure, the representability of local Shimura varieties themselves remains an open question outside the Abelian type. Furthermore, the conjecture that these local Shimura varieties are perfectoid at the infinite level has not been proven for any case beyond the Abelian type.\\

In this paper, we bridge this gap by shifting the perspective and proving that certain exceptional Shimura varieties can be embedded into Siegel modular varieties at the level of the varieties themselves, even when a group-theoretic embedding of the data is impossible. This ``variety-level'' embedding allows us to resolve the following conjectures on representability and perfectoidness for an infinite family of non-Abelian type Shimura data:

\begin{conjecture*}[\hypertarget{representintegrallocal}{A}]
    Let $\left(G, \mu, b\right)$ be an integral local Shimura datum. Then the integral model for the local Shimura variety at hyperspecial level, $\mathcal{M}^{\textup{int}}_{(G, \mu, b)}$ , is representable by a smooth formal scheme.
\end{conjecture*}

\begin{conjecture*}[\hypertarget{perfectoidintegrallocal}{B}]
    Let $\left(G, \mu, b\right)$ be a local Shimura datum. Then the local Shimura variety at infinite level, $\mathcal{M}_{(G, \mu, b),\infty} = \varprojlim \mathcal{M}_{(G,\mu, b),K}$, is representable by a perfectoid space.
\end{conjecture*}

\begin{conjecture*}[\hypertarget{perfectoidglobal}{C}]
    Let $(G, \mu)$ be a Shimura datum, and let $K^p \subset G(\mathbb{A}^{\infty, p})$ be a fixed level away from $p$. The inverse limit
    \[ \varprojlim_{K_p} \Sh_{K_p K^p}(G, \mu) \]
    over compact open subgroups $K_p \subset G(\mathbb{Q}_p)$, corresponding to the levels at $p$, is representable by a perfectoid space.
\end{conjecture*}

Our work was motivated by Bültel's construction of integral models of a class of Shimura varieties, called ``meta-unitary'' Shimura varieties, via the theory of displays \cite{bultel2008pel}. Note that the class of meta-unitary Shimura varieties includes those for some exceptional groups; these Shimura varieties are not of Abelian type. In particular, we were interested in proving Langlands-Rapoport conjecture and Conjecture \hyperlink{perfectoidglobal}{C}. We have also defined the local counterparts of the meta-unitary Shimura varieties and proved Conjectures \hyperlink{representintegrallocal}{A} and \hyperlink{perfectoidintegrallocal}{B} for them.\\

As noted above, the main difficulty in working with Shimura varieties that are not of Abelian type is the lack of a moduli description. Let $(G,\mu)$ be a Shimura datum that is not of Abelian type. When we embed $G$ in some $\GL_n$, the cocharacter $\mu$ does not remain minuscule; therefore, in the global setting, the associated Hodge structures do not correspond to Abelian varieties. In the local setting, the associated crystals are not related to $p$-divisible groups.\\

The idea of B\"ultel's construction is to define a ``Dieudonn\'efication'' functor called \emph{$\Flex$}, which maps certain graded non-minuscule crystals to $p$-divisible groups and uses it to describe the special fiber of some Shimura varieties. He also defines a functor from some graded non-minuscule Hodge structures to Hodge structures with weights $(1,0)$ and $(0,1)$. These two functors help to relate the mod $p$ points of $\Sh(G,\mu)$ to Abelian varieties and obtain a map from the special fiber of $\Sh(G,\mu)$ to a Siegel modular variety.\\

The main innovation of this paper, inspired by prismatic geometry and mixed-characteristic local shtukas, is the construction of a mixed-characteristic functor that specializes to both functors defined by B\"ultel. Conceptually, we use the framework of prismatic crystals and shtukas with several legs to decompose non-minuscule cocharacters into a sequence of minuscule ones. Then, by interpreting shtukas as modifications of vector bundles on the Fargues-Fontaine curve, we combine these legs to obtain a shtuka with a single leg. For a suitable class of Shimura data —referred to as \emph{regularly sparse Shimura data}— the resulting shtuka is bounded by a minuscule cocharacter, and no information is lost in the process. This helps us embed these Shimura varieties into the Siegel modular varieties and embed their local counterparts inside Rapoport-Zink spaces. Using this approach, we prove some special cases of the above conjectures.\\

We prove the above conjectures for a large class of Shimura varieties, including the following: those associated with data $(G, \mu)$ splitting over $\mathbb{Q}_{p^f}$, where $G$ is of one of the types $B_n,C_n$, or $E_7$; the derived group $G^{\der}$ is simply connected; the adjoint quotient $G^{\ad}$ is simple and quasi-split over $\mathbb{Q}_p$; and $G^{\ad}\otimes\mathbb{R}$ has exactly one non-compact factor and at least five compact factors.\\

We now provide further technical details on the definition of the ``Dieudonn\'efication'' functor $\Flex$. To illustrate the main ideas and challenges involved, we consider a simple special case. Suppose $M = \oplus_{i\in \frac{\mathbb{Z}}{r}} M_i$ is a graded finitely generated projective module over the ring of $p$-typical Witt vectors, $W(k)$, where $k$ is a perfect field of characteristic $p$. Let $F\colon M \to M$ be a degree one Frobenius-linear map. Assume that 
$F_i\coloneqq F_{\mid_{M_i}}\colon M_i \to M_{i+1}$ is an isomorphism for $i \neq 0$, and $p^w M_{1} \subset F_0 M_0$ for some $w \leq r$. Consider the graded module $\tilde{M} = \oplus_{i\in \frac{\mathbb{Z}}{r}} \tilde{M}_i$, where $\tilde{M}_0 = M_0$ and for $i \not= 0$ $$\tilde{M}_{i} = F^{i} M_0 + p^{\max\{w-i, 0\}} M_{i}$$ It is straightforward to verify that $p\tilde{M} \subset F\tilde{M}$, and hence $\tilde{M}$ naturally acquires the structure of a Dieudonn\'e module. In fact, the ``empty'' components $M_i$ for $i \neq 0$ are used to decompose the non-minuscule cocharacter —encoding the relative position of the lattices $FM_0$ and $M_1$— into a collection of minuscule cocharacters. To extend this construction over a general base of characteristic $p$, B\"ultel lifted it to the framework of displays.\\

Let us remark on why these graded crystals appear naturally. Assume that $M$ is a crystal over $\mathbb{F}_p$ equipped with an action of $W(\mathbb{F}_{p^r})$. Write $W(\mathbb{F}_{p^r}) = W(\mathbb{F}_p)[\alpha]$ (note that the actions of $\alpha$ and $F$ commute with each other). Then, one can prove that $M_{\mathbb{F}_{p^r}}$ decomposes into eigenspaces corresponding to the action of $\alpha$. This decomposition induces a natural $\mathbb{Z}/r\mathbb{Z}$-grading on $M$, with the Frobenius morphism having degree one. Conversely, by descent, from a graded crystal over $\mathbb{F}_{p^r}$, one obtains a crystal over $\mathbb{F}_p$ endowed with an action of $W(\mathbb{F}_{p^r})$. This is the motivation for working with graded crystals.\\

The naive way to extend this construction to mixed-characteristic is to replace $p$ in the above formula by $d: = V\left(1\right)$, where $V$ is the Verschiebung of the Witt ring. However, in general, $V\left(1\right)$ is not Frobenius-invariant, so we only obtain a map
$$F^iM_0+d^{\max\{w-i,0\}}M_{i}\to F^{i+1}M_0+\sigma\left(d\right)^{\max\{w-i-1,0\}}M_{i+1}$$
which is an isomorphism after inverting $\sigma\left(d\right)$ (here $\sigma$ is the Frobenius of the Witt ring). We obtain a Frobenius crystal with several legs, and here is the challenge.\\

More generally, one can make this construction for prismatic $F$-crystals. In the perfectoid world, we set $$\tilde{M}_i=F^iM_0+\sigma^i(d)^{\max\{w-i, 0\}}M_i$$ This yields a shtuka with legs at $d, \sigma\left(d\right), \dots, \sigma^{w-1}\left(d\right)$. However, all these points share the same untilt and map to a single point in the Fargues-Fontaine curve. Thus, one expects that these legs can be combined into a single leg, and indeed, this is possible. If $w$ is even, when we combine the legs, we obtain a shtuka with one leg bounded by a minuscule cocharacter.\\

We define the stack of graded prismatic displays with several legs and construct the functor \emph{$\Flex$} in two steps. First, we restrict our attention to graded prismatic displays for which the Frobenius map is an isomorphism in all but a few grades. These ``empty'' components are essential. The Frobenius map at some grades is bounded by some non-minuscule cocharacter, and we will use the empty spaces to break it into Frobenius maps in several grades, each bounded by a minuscule cocharacter, such that the sum of these minuscule cocharacters will be the original non-minuscule cocharacter. To isolate such displays, we begin by defining an auxiliary stack (see Definition \ref{auxiliarystackdef}). Then, using the constructions introduced in the preceding paragraphs, we define a functor from this auxiliary stack to the stack of graded prismatic displays with several legs. This functor lifts B\"ultel’s construction from characteristic $p$ (the special fiber) to mixed characteristic.\\

Next, considering shtukas as modifications of vector bundles on the Fargues--Fontaine curve, we combine these legs into a shtuka with a single leg. We show that for regularly sparse Shimura data, the resulting shtuka is bounded by a minuscule cocharacter.\\

We prove that the functor $\Flex$ is an isogeny modulo $p$, which allows us to lift the functor $\Flex$ to the local Shimura variety. One also observes that, over the generic fiber, $\Flex$ is an isomorphism near the closed subspace of points with radius $0$ of the curve $\mathcal{Y}$ (this is the curve whose quotient by Frobenius yields the Fargues-Fontaine curve). This allows us to lift the functor $\Flex$ to the moduli space of shtukas with deeper levels over the generic fiber.\\

We also define the class of \emph{sparse local Shimura varieties}. Many local Shimura varieties in this class are not of local Hodge type; in fact, the class includes local Shimura varieties associated with exceptional groups, such as those of type $E_7$. Using the construction described above, we obtain a map from the integral model of regularly sparse local Shimura varieties with hyperspecial level to Rapoport--Zink spaces. By applying the representability of the formal completion of the integral local Shimura variety at closed points, along with a criterion due to Pappas and Rapoport, we prove that this map is representable by a closed immersion of formal schemes.\\

\begin{thm}\label{regularsparseisrepresentable}
    Let $\left(G,\mu,b\right)$ be a regularly sparse Shimura datum. The associated integral local Shimura variety $\mathcal{M}^{\textup{int}}_{(G,\mu,b)}$, considered as a v-sheaf, is representable by the diamond of a smooth formal scheme.
\end{thm}

This result allows us to establish the representability of integral local Shimura varieties with hyperspecial level in several new cases.

\begin{thm}\label{regularsparseperfect}
    Let $\left(G,\mu,b\right)$ be a regularly sparse Shimura datum. The local Shimura variety at infinite level $\mathcal{M}_{(G,\mu,b),\infty}=\varprojlim_K \mathcal{M}_{(G,\mu,b),K}$ is representable by a perfectoid space.
\end{thm}

Furthermore, over the generic fiber, we can lift the functor $\Flex$ and obtain a closed immersion of the local Shimura varieties for regularly sparse Shimura datum to the Rapoport-Zink space with deeper levels. This provides a proof of Conjecture \hyperlink{perfectoidintegrallocal}{B} in some new cases.\\

For a global application, we use the functor $\Flex$ defined above to construct a closed immersion from regular meta-unitary Shimura varieties (see Appendix \ref{metaunitarydatumdfn} for the definition) to Siegel modular varieties. More precisely, we have the following:

\begin{cor}\label{closed immersion}
    Let $(G,\mu)$ be a regular meta-unitary Shimura datum, and $K^p\subset G(\BA^{p,\infty})$ a fixed level away from $p$. Then, there exist an open compact subgroup $K\subset G(\BQ_p)$, a Siegel Shimura datum $(G^\triangleright,\mu^\triangleright)$ (explicitly constructed), a level $\tilde{K}\subset G^\triangleright(\BA^\infty)$, and a closed immersion $$\Flex\colon\Sh_{K_pK^p}\left(G,\mu\right)\xhookrightarrow{\quad} \Sh_{\tilde{K}}\left(G^\triangleright,\mu^\triangleright\right).$$
\end{cor}

We conjecture that the level $K_pK^p$ is hyperspecial.\\

Write $K=K_pK^p$ where $K_p$ is the level at $p$ and $K^p$ is the level away from $p$. If we let the level at $p$ vary, we obtain a pro-system of \'etale covers of $\Sh_K(G,\mu)$. This tower gives us a $G(\mathbb{Z}_p)$-torsor over $\Sh_K(G,\mu)$. Hansen \cite{Hansen2016PeriodMA} proved that this $G(\mathbb{Z}_p)$-torsor corresponds to a shtuka on the Shimura variety. We use the existence of this shtuka over Shimura varieties with deeper levels to lift the above closed immersion to deeper levels. Under the technical assumption (CS) that the display defined by B\"ultel on the meta-unitary Shimura varieties is the same as the one defined by Hansen, we prove Conjecture \hyperlink{perfectoidglobal}{C} for these Shimura varieties. Note that we expect this assumption to be true, and we will return to proving it in our future work \cite{hedayatzadehrappaportlanglands}.

\begin{thm}\label{globalperfectoid}
    Let $\left(G,\mu\right)$ be a regular meta-unitary Shimura datum and $K^p\subset G(\BA_{\BQ}^{p,\infty})$ be a compact open subgroup. Under the assumption (CS), the Shimura variety at infinite level $$\Sh_{K^p}\left(G,\mu\right)\coloneqq\varprojlim_{K_p}\, \Sh_{K_pK^p}\left(G,\mu\right)$$ is representable by a perfectoid space.
\end{thm}

These theorems are just some applications of the existence of the embedding of meta-unitary Shimura varieties into the Siegel space. Indeed, one can naturally expect that by using this embedding, one can generalize many results from the case of Shimura varieties of Abelian type to meta-unitary Shimura varieties.

We would also like to address the following two questions regarding this work:

\begin{que}
    Can one lift this closed immersion to parahoric levels on the integral model?
\end{que}

\begin{que}
    Is it possible to use these results to prove the representability of the integral model of a more general class of Shimura varieties?
\end{que}

\subsection*{Structure of the Paper}

In Section \ref{gradeddisplay}, we define the stack of graded prismatic displays with multiple legs and a twisted version of this stack. These stacks lift the stacks defined by B\"ultel in \cite{bultel2008pel} and help us define the functor $\Flex$. The definitions related to prismatic displays are reviewed in Appendix \ref{PrismDispGrp}. In Section \ref{definieflex}, we define the functor $\Flex$ over graded shtukas. In Section \ref{sectionsparse}, we define the class of \emph{sparse Shimura data} and a ``non-standard'' embedding from the associated integral local Shimura variety to the stack of shtukas with one leg, bounded by some cocharacter. If the Shimura datum satisfies some extra conditions, the cocharacter will be minuscule. A local Shimura datum satisfying these extra conditions is said to be \emph{regularly sparse}. In Section \ref{regularlysparsesection}, we prove that the integral model of regularly sparse Shimura varieties is representable and that it is perfectoid at the infinite level. In Section \ref{integralmodelsection}, we construct a closed embedding from the integral model of regular meta-unitary Shimura varieties into Siegel varieties. In Appendix \ref{defmeta-unitarysection}, we review the definition of meta-unitary Shimura varieties and the construction of their integral model. In Section \ref{deeperlevelsection}, we prove that regular meta-unitary Shimura varieties at infinite level are perfectoid. The appendices provide necessary background on prismatic display groups, meta-unitary Shimura varieties, shtukas, and local Shimura varieties.

\subsection*{Acknowledgements}
We are grateful to Oliver B\"ultel, Miaofen Chen, and Aliakbar Hosseini for many helpful discussions. We thank Kiran Kedlaya for his interest and for answering some of our questions.

\subsection*{Notations and Conventions}\label{Not&Conv}

Here is a list of notations and conventions that we will use throughout the paper. 

\begin{itemize}
    \renewcommand\labelitemi{--}
    \item $p$ is a prime number and $q$ is a power of $p$.
    \item Unless specified to the contrary, all rings are commutative with a $1$.
    \item For a ring $R$, we denote by $W\left(R\right)$ the ring of ($p$-typical) Witt vectors, and by $W_n(R)$ the truncated Witt vectors of length $n+1$.
    \item $\BZ_{p^f}\coloneqq W\left(\BF_{p^f}\right)$ and $\BQ_{p^f}\coloneqq\BZ_{p^f}[1/p]$.
    \item We will use $\CG$ for a reductive group over some unramified extension of $\BZ_p$ (which will be specified in the context) and $G$ for a reductive group over some field extension of $\BQ$. Typically, $G$ will be the generic fiber of $\CG$.
    \item For a ring $R$ of characteristic $p$, we denote by $A_{\textup{crys}}(R)$ the crystalline period ring of $R$. 
    \item For a $p$-adic ring $A$, $\sigma_A$ is a Frobenius lift from $A/p$ to $A$.
    \item Let $R$ be a ring and $S$ an $R$-algebra. The Weil restriction functor from $S$ to $R$ is denoted by $\Res_{S/R}$.
    \item $\Perf$ is the site of perfectoid spaces in characteristic $p$ with the \'etale topology.
    \item $\mathbb{A}$ is the ring of adeles of $\mathbb{Q}$, $\mathbb{A}^\infty$ is the ring of finite adeles, and $\mathbb{A}^{\infty,p}$ is the ring of finite adeles away from $p$.
    \item For an affinoid perfectoid space $\Spa(S,S^+)$, we denote by $\mathcal{Y}_S$ the curve $\Spa W(S^+)/([w],p)$ and $\kappa\colon\mathcal{Y}_S\to \mathbb{R}$ is the radius map.
    \item For an interval $I\subset\BR$, we denote by $\mathcal{Y}_{S,I}$ the inverse image $I$ under $\kappa$.
    \item If $\mathcal{P}$ is a vector bundle or torsor on $\mathcal{Y}_S$, we will denote the restriction of $\mathcal{P}$ to $\mathcal{Y}_{S,I}$ by $\mathcal{P}_{I}$.
    
\end{itemize}

\normalsize

\section{Graded Prismatic Displays}\label{gradeddisplay}

For the applications we have in mind, we need to introduce some graded variants of the stack of prismatic displays defined in \cite{hedayatzadeh2024deformations}. The first variant allows us to have different ``legs'' at different grades. Let $r$ be a natural number.

\begin{dfn}
    A \emph{multi-leg prism} $(A,\mathbf{I}
)$ is a $\delta$-ring $A$ together with a vector $\mathbf{I}\coloneqq\left(I_i\right)_{i
    \in \BZ/r}$, with $I_i$ an ideal of $A$, such that for each $i\in \BZ/r$, the pair $\left(A,I_i\right)$ is an (oriented) prism.
\xqed{\blacktriangle}\end{dfn}

\begin{dfn}
    A \emph{graded Breuil-Kisin-Fargues module} over a multi-leg prism $(A,\mathbf{I})$ is a pair $\left(M,F\right)$, where $M=\bigoplus_{i\in\BZ/r}M_i$ is a graded finitely generated projective $A$-module and $F\colon\sigma^*M\to M$ is an $A$-linear morphism of degree $1$ (mod $r$) such that for all $i\in\BZ/r$ the base change of $$F_i\coloneqq F_{\mid \sigma^*M_i}\colon\sigma^* M_i\to M_{i+1}$$ to $A[\frac{1}{I_{i}}]$ is an isomorphism. We call a graded Breuil-Kisin-Fargues module \emph{banal} if all $M_i$ are free $A$-modules.

    Let $\CG$ be a reductive group over some unramified extension of $\BZ_p$. A \emph{graded Breuil-Kisin-Fargues module with $\CG$-structure} (\emph{$\CG$-BKFg module} for short) over $(A,\mathbf{d})$ is a pair $\left(\SM,F\right)$, where $\SM=\bigoplus_{i\in\BZ/r}\SM_i$ and $F=\oplus_{i\in\BZ/r} F_i$, and for each $i\in\BZ/r$, $\SM_i$ is a $\CG$-torsor over $A$ and $F_i\colon\sigma^*\SM_{i}\to \SM_{i+1}$ is a $\CG$-equivariant morphism over $A$ such that the base change of $F_i$ to $A[\frac{1}{d_{i}}]$ is an isomorphism. We call a $\CG$-BKFg module \emph{banal} if all $\SM_i$ are the trivial $\CG$-torsor. A morphism from $(\SM,F)$ to $(\SM',F')$ is a morphism of graded $\CG$-torsors that commutes with $F$.
\xqed{\blacktriangle}\end{dfn}

\begin{example}
    Consider the multi-leg prism $\left(A,\left((d),(\sigma(d))\right)\right)$ and set $M=A^{\oplus 2}\oplus A^{\oplus 2}$.
    Let $F=F_1\oplus F_2$ be the Frobenius map of degree $1$ defined by
    $$F_1=\begin{pmatrix}
    d & 0\\
    0 & 1
    \end{pmatrix}\quad\text{and}\quad F_2=\begin{pmatrix}
    \sigma(d) & 0\\
    0 & 1
    \end{pmatrix}$$
    Then $(M,F)$ is a graded Breuil-Kisin-Fargues module over $\left(A,\left((d),(\sigma(d))\right)\right)$. 
\xqed{\blacksquare}\end{example}

\begin{rmk}
    There is a Tannakian description of $\CG$-BKFg modules: this category is equivalent to the category of fiber functors from the category of representations of $\CG$ to the category of graded Breuil-Kisin-Fargues modules.
\xqed{\lozenge}\end{rmk}

\begin{dfn}
    Assume that $\left(A,\mathbf{I}\right)$ is a multi-leg prism. Set $${}_+\vec{\CG}\left(A,\mathbf{I}\right)\coloneqq\prod_{i\in\BZ/r} \CG\left(A\right)$$ and the \emph{loop group} $$\vec{\CG}\left(A,\mathbf{I}\right)\coloneqq\prod_{i\in\BZ/r} \CG\left(A\left[\frac{1}{I_{i}}\right]\right)$$
\xqed{\blacktriangle}\end{dfn}

\begin{rmk}
    Note that after fixing a trivialization, a banal $G$-BKFg module is given by an element of $LG$. Therefore, the groupoid of banal $G$-BKFg modules is the quotient groupoid $\left[\vec{\CG}\left(A,\mathbf{I}\right):{}_+\vec{\CG}\left(A,\mathbf{I}\right)\right]$ where $k\in {}_+\vec{\CG}\left(A,\mathbf{I}\right)$ acts on $\vec{\CG}\left(A,\mathbf{I}\right)$ by $$\left(k_i\right)\cdot\left(U_i\right)\coloneqq\left(k_i^{-1}U_i\sigma\left(k_{i-1}\right)\right)$$
\xqed{\lozenge}\end{rmk}

We want to introduce the stack of graded displays, which has a natural functor to the stack of graded Breuil-Kisin-Fargues modules.

\begin{dfn}
    Let $\left(A,\mathbf{I}\right)$ be a multi-leg prism. The \emph{relative prismatic site} over $\left(A,\mathbf{I}\right)$ is the site of multi-leg prisms $\left(B,\mathbf{I}\right)$ with completely \'etale topology. We denote this site by $\left(A,\mathbf{I}\right)_\prism$.
\xqed{\blacktriangle}\end{dfn}

    Let $\left(A,\mathbf{I}\right)$ be a multi-leg prism, $\CG$  a reductive group over some unramified extension of $\mathbb{Z}_p$, and $\boldsymbol{\mu}=(\mu_{i})$ be a cocharacter of $\prod_{i\in\BZ/r} \CG$ defined over that extension. In Appendix \ref{PrismDispGrp}, we define the associated canonical limit group $\CG(A,\mathbf{I})$, the prismatic display group $\CG^{\boldsymbol{\mu}}(A,\mathbf{I})$, and an action of $\CG^{\boldsymbol{\mu}}(A,\mathbf{I})$ on $\CG(A,\mathbf{I})$, called the \emph{$\Phi_{\mathbf{I}}^{\boldsymbol{\mu}}$-conjugation}.

\begin{dfn}
    Let $\left(A,\mathbf{I}\right)$ be a multi-leg prism. The groupoid of \emph{banal graded prismatic higher $\left(\CG,\boldsymbol{\mu}\right)$-displays} over $\left(A,\mathbf{I}\right)$ is the quotient groupoid $$\mathcal{B}_{\left(\CG,\boldsymbol{\mu}\right)}^{\textup{ban}}\left(A,\mathbf{I}\right)\coloneqq\left[\CG\left(A,\mathbf{I}\right):\CG^{\boldsymbol{\mu}}\left(A,\mathbf{I}\right)\right]$$ where $\CG^{\boldsymbol{\mu}}\left(A,\mathbf{I}\right)$ acts on $\CG\left(A,\mathbf{I}\right)$ by $\Phi_{\mathbf{I}}^{\boldsymbol{\mu}}$-conjugation. 
\xqed{\blacktriangle}\end{dfn}

 \begin{dfn}\label{DefPrimStackPrimDisp}
    Define a presheaf of groupoids by sending a multi-leg prism $\left(B,\mathbf{I}\right)$ over $\left(A,\mathbf{I}\right)$ to $\mathcal{B}_{\left(\CG,\boldsymbol{\mu}\right)}^{\textup{ban}}\left(C,\mathbf{I}\right)$. This presheaf is, in fact, a prestack (because both objects and morphisms are sheaves).  We define the stack of \emph{graded prismatic higher $\left(\CG,\boldsymbol{\mu}\right)$-displays}, or \emph{gph $(\CG,\boldsymbol{\mu})$-displays} to be the sheafification of this presheaf and denote it by $\mathcal{B}^{\prism_{/\left(A,\mathbf{I}\right)}}_{\left(\CG,\boldsymbol{\mu}\right)}$.
\xqed{\blacktriangle}\end{dfn}

\begin{pro}\label{fiberfargueskisinmodule}
    Assume that $\left(A,\mathbf{I}\right)$ is a multi-leg prism. There is a fully-faithful functor from $\mathcal{B}^{\prism_{/\left(A,\mathbf{I}\right)}}_{\left(\CG,\boldsymbol{\mu}\right)}$ to the stack of $\CG$-BKFg modules over $\left(A,\mathbf{I}\right)$
\end{pro}

\begin{proof}
    It is enough to define the functor over banal objects. To do this, we send $\left(U_i\right)\in \CG\left(A,\mathbf{I}\right)$ to $\left(U_i\mu_{i-1}\left(I_{i-1}\right)\right)\in \vec{\CG}\left(A,\mathbf{I}\right)$ and a morphism $k$ to $k$. One can check that this functor is well defined and fully faithful, similar to \cite[Proposition 4.11]{hedayatzadeh2024deformations}.
\end{proof}

Let $\CG_0$ be a reductive group scheme over $\mathbb{Z}_{p^r}$ and set $\CG=\operatorname{Res}_{\mathbb{Z}_{p^r}/\mathbb{Z}_p}\CG_0$. There is a canonical isomorphism $\CG_{\mathbb{Z}_{p^r}} \cong \prod_{\BZ/r} \CG_0$.
 Let $\boldsymbol{\mu}\colon \mathbb{G}_m \to \CG_{\mathbb{Z}_{p^r}}$ be a cocharacter. Under the above isomorphism, $\boldsymbol{\mu}$ can be written as a cocharacter $\boldsymbol{\mu}=(\mu_{i})$ for $\prod_{\BZ/r} \CG_0$.

\begin{pro}
    Assume that $I_i=I$ for all $i$ and $(A,\mathbf{I})$ is a $\BZ_{p^r}$-prism, then $\mathcal{B}^{\prism_{/\left(A,\mathbf{I}\right)}}_{\left(\CG_0,\boldsymbol{\mu}\right)}$ and $\mathcal{B}^{\prism_{/\left(A,I\right)}}_{\left(\CG,\boldsymbol{\mu}\right)}$ are isomorphic.
\end{pro}

\begin{proof}
   Note that there is an isomorphism $A\otimes \mathbb{Z}_{p^r}\xrightarrow{\sim}\prod_{\BZ/r}A$, which sends $a\otimes \alpha$ to $\left(a\sigma^i(\alpha\right))$, where $\sigma$ is the Frobenius automorphism of $\BZ_{p^r}$. Under this isomorphism, $\sigma_A(a)\otimes \alpha$ maps to $\left(\sigma_A(a)\sigma^{i}(\alpha)\right)$. Therefore, under the isomorphism of group schemes $\iota\colon\CG_{\mathbb{Z}_{p^r}}\xrightarrow{\sim}\prod_{\BZ/r} \CG_0$, if $g\in \CG(A)$ and $\iota\left(g\right)=(g_i)_{i\in\BZ/r}$, we have $\iota\left(\sigma_A\left(g\right)\right)=\left(\sigma_A(g_{i-1})\right)_{i\in\BZ/r}$. It follows that $\iota$ is $\CG^{\boldsymbol{\mu},r}$-equivariant and induces an isomorphism between the quotient groupoids.
\end{proof}

Let $\Div^1$ be the v-sheaf that sends an affinoid perfectoid space $S=\Spa(R,R^+)$ to the set of untilts of $S$, in other words, the set of Cartier divisors on $\mathcal{Y}_S$, and denote $(\Div^1)^{\times r}$ by $\Div^r$. So, a map $\Spa\left(R, R^+\right) \to \Div^r$ is the data of $r$ untilts $(R^\sharp_i,R_i^{\sharp +})$ of $(R,R^+)$. Let $I_i=\ker\left( \theta_i\colon W(R^{+}) \to R_{i}^{\sharp+}\right)$. Note that there is a natural bijection between the points of $\Div^1(R)$ and distinguished ideals $I\subset W(R)$

\begin{dfn}\label{DefPerfDisp}
    Over the site $\Perf /\Div^r$, consider a presheaf of groupoids that sends an affinoid $\Spa\left(R,R^+\right)\to \Div^r$, to $$\mathcal{B}_{\left(\CG_0,\boldsymbol{\mu}\right)}^{\textup{ban}}\left(W(R^+),\left(I_i\right)_{i=0}^r\right)$$ We denote the sheafification of this presheaf by $\CB^{\textup{perf},r}_{\left(\CG_0,\boldsymbol{\mu}\right)}$.
\xqed{\blacktriangle}\end{dfn}

\begin{rmk}\label{compprismdisplay}
    By a similar argument to \cite{ito2023prismatic}, one can show that the natural map between $\CB^{\textup{perf},r}_{\left(\CG_0,\boldsymbol{\mu}\right)}(\left(\Spa(R,R^+)\to\Div^r\right))$ and $\CB^{\prism}_{\left(\CG_0,\boldsymbol{\mu}\right)}(A_{\textup{inf}}(R),\mathbf{I})$ is an isomrphism.
\xqed{\lozenge}\end{rmk}

Let us now provide an alternative definition of graded prismatic higher $\left(\CG_0,\boldsymbol\mu\right)$-displays, using torsors:

\begin{pro}\label{PropDispTors}
    Let $S$ be an object of $\Perf/\Div^r$. Giving a gph $\left(\CG_0,\boldsymbol\mu\right)$-display over $S$ is the same as giving a pair $\left(\SQ,\alpha\right)$, where $\SQ$ is a torsor for the sheaf $\CG_0^{\boldsymbol{\mu}}$ over $S$ and $\alpha\colon\SQ\to {}_+\vec{\CG}_0$ is a $\CG_0^{\boldsymbol{\mu}}$-equivariant morphism.
\end{pro}

\begin{proof}
    The proof is similar to the proof of \cite[Proposition 4.20]{hedayatzadeh2024deformations}, by adapting it to the multi-leg setting.
\end{proof}

\begin{rmk}
    A graded $\CG$-shtuka over $S$ is a graded sheaf $\oplus \mathcal{P}_i$ over $\mathcal{Y}_S$ such that each $\mathcal{P}_i$ is a $\CG$-torsor together with Frobenius maps $F_i\colon\sigma^*\mathcal{P}_i\to \mathcal{P}_{i-1}$ that are isomorphism away from legs.
\xqed{\lozenge}\end{rmk}

\begin{rmk}
    We say that a $\CG$-shtuka $\mathcal{P}$ over $S$ is bounded by $\mu$ at leg $y$, if for any geometric point of $S$ like $C$, after choosing a trivialization of $\mathcal{P}$ near $y$, the Frobenius is represented by an element of $\CG(B_{dR}^+(y))\mu(d)\CG(B_{dR}^+(y))$.
\xqed{\lozenge}\end{rmk}

\begin{construction}\label{ConsFib}
    Let us denote by $\Sht^r_{(\CG, \boldsymbol{\mu})}$ the stack of $\CG$-shtukas with $r$ legs. We want to define a functor
    \[
    \textup{Fib}_{\CG} \colon \mathcal{B}^{\textup{perf},r}_{(\CG_0, \boldsymbol{\mu})} \to \Sht^r_{(\CG, \boldsymbol{\mu})}
    \]
    By Proposition \ref{fiberfargueskisinmodule}, we have a functor from $\mathcal{B}^{\textup{perf},r}_{\left(\CG_0, \boldsymbol\mu\right)}$ to the category of $\CG$-BKFg modules, and so it is enough to define a functor from the latter category to $\Sht^r_{(\CG, \boldsymbol{\mu})}$. It suffices to define this functor on banal objects. We send $\left(U_i\right)$ to the graded $\CG_0$-shtuka defined by the trivial graded $\CG_0$-torsor and the Frobenius morphism $\left(U_i\right)$.
\xqed{\blacktriangledown}\end{construction}

We also want to define another variant that will be useful in the next section. To do this, we first need to introduce some combinatorial notions. To motivate them, let us recall the ``Dieudonnéfication'' functor from the introduction via a simple example:

\begin{example}\label{motivatingexm}
    Let $k$ be a perfect field in characteristic $p$. Suppose we have a graded crystal \( M = \bigoplus_{i\in\BZ/r} M_i \) over $k$, with a degree one Frobenius $F=\oplus F_{i}\colon\sigma^*M\to M$ such that
    \[
    p^w M_1 \subset F_0 \left(\sigma^* M_0\right) \subset M_1,
    \]
    for some $w\in\BN$ and \( F_i \) is an isomorphism for all \( i \neq 0 \). Our goal is to convert this into a graded crystal \(\tilde{M} = \bigoplus_{i\in\BZ/r} \tilde{M}_i\) satisfying
    \[
    p \tilde{M}_1 \subset F_0 (\sigma^*\tilde{M}_0) \subset \tilde{M}_1
    \] so as to have a (graded) Dieudonné module.
    To achieve this, we set
    \[
    \tilde{M}_i \coloneqq F^i (\sigma^*M_0) + p^{\max\{w - i, 0\}} M_i.
    \]
    
    In the mixed characteristic case, a similar construction applies. Let $(A,d)$ be an oriented prism and let \( M = \bigoplus_{i\in\BZ/r} M_i \) be a graded projective module over $A$ together with a degree one Frobenius $F=\oplus F_{i}\colon\sigma_A^*M\to M$ such that 
    \[
    d^w M_1 \subset F_0 (\sigma_A^* M_0) \subset M_1,
    \]
    for some $w\in\BN$ and \( F_i \) is an isomorphism for all \( i \neq 0 \). We set 
    \[
    \tilde{M}_i \coloneqq F^i (\sigma^*M_0) + \sigma^{i-1}(d)^{\max\{w - i, 0\}} M_i
    \] and $\tilde{M}\coloneqq\bigoplus_{i\in\BZ/r}\tilde{M}_i$.
    Then $\tilde{M}$ is a graded Breuil-Kisin-Fargues module over the multi-leg prism $(A,\mathbf{d})$, where $\mathbf{d}=(d,\sigma(d),\dots,{\sigma^{w-1}}(d),d,\dots,d)$, and is bounded by a minuscule cocharacter. Note that $\tilde{M}$ has legs at $d,\sigma(d),\dots,{\sigma^{w-1}}(d)$.
\xqed{\blacksquare}\end{example}

To make this construction in the language of displays and to handle more complicated examples, we need to distinguish between the grades where Frobenius is an isomorphism and those where it is not. Let \(\Gamma\) denote the set of grades where Frobenius is \emph{not} an isomorphism.

The grades on which Frobenius \emph{is} an isomorphism serve as our ``empty space'', allowing us to decompose the non-minuscule cocharacters that encode the relative position of \(M_\gamma\) and \(F (\sigma^*M_\gamma)\) for \(\gamma \in \Gamma\). For this, we assign some of the empty grades to each \(\gamma \in \Gamma\), and we also keep track of the distance between these empty grades and \(\gamma\).

\begin{setup}\label{notationgamma}
    Consider a nonempty subset $\Gamma\subset \mathbb{Z}/r\mathbb{Z}$. Set $\CG^\Gamma= \prod_{\gamma\in \Gamma} \CG_0$ and let $\mu^\Gamma=(\mu^{\gamma})_{\gamma\in\Gamma}$ be a cocharacter of $\CG^\Gamma$. Consider the natural bijection $$\textup{mod}\colon\mathbb{Z}/r\mathbb{Z}\xrightarrow{\sim}\{1,2,\cdots,r\}$$ and let $\gamma_i$ be the $i$-th element of $\Gamma$ with the ordering induced from this bijection. Set $$r^{+}_i\coloneqq\textup{mod}\left(\gamma_{i+1}-\gamma_{i}\right)\in\BN$$ For $1\le i\le r$, let $[\gamma_i,\gamma_{i+1})\subset \BZ/r\BZ$ be the subset of elements $w\in\BZ/r\BZ$, such that $\gamma_i\le w<\gamma_{i+1}$ if $i\le |\Gamma|$ and $\gamma_i<w<$ or $w<\gamma_{1}$ if $i=|\Gamma|$ (i.e., elements that are strictly between $\gamma_i$ and $\gamma_{i+1}$ mod $r$). For $w\in [\gamma_i,\gamma_{i+1})$, define $t^+\left(w\right)\coloneqq\textup{mod}\left(w-\gamma_i\right)$.
\end{setup}

It is convenient to first define an auxiliary stack that allows us to focus only on the grades where Frobenius is \emph{not} an isomorphism. With the notations of Example \ref{motivatingexm}, the idea is to compose morphisms \( (\sigma^{r_i^{+}-t^+(w)-1})^*F_w \) for \( w \in [\gamma_i, \gamma_{i+1}) \) (note that all of these morphisms are isomorphisms except the first one) to obtain a morphism
\[
(\sigma^{r_{i}})^*M_{\gamma_i} \to M_{\gamma_{i+1}}.
\]

To do this in a group theoretic way, we give the following definition:

\begin{dfn}\label{auxiliarystackdef}
    Define the prestack of groupoids sending the prism $\left(A, I\right)$ to the quotient groupoid $$\left[\CG^\Gamma\left(A\right) / \prod_{\gamma \in \Gamma} \CG^{\mu_{\gamma}}\left(A, I\right)\right]$$ where $\prod_{\gamma \in \Gamma} \CG^{\mu_{\gamma}}\left(A, I\right)$ acts on $\CG^\Gamma\left(A\right)$ by the formula: $$\left(k_{\gamma_i}\right)_{\gamma_i\in \Gamma}.\left(U_{\gamma_i}\right)_{\gamma_i\in\Gamma}=\left( k_{\gamma_i}^{-1}U_{\gamma_i} \Phi_{\mu_{i-1}}^{\sigma^{r^+_{i-1}}(I)}\left(\sigma^{r^+_{i-1}-1}\left(k_{\gamma_{i-1}}\right)\right)\right)_{\gamma_i\in \Gamma}$$

    We define the stack \emph{$\CB^{\textup{perf}}_{\left(\CG^\Gamma,\mu^\Gamma\right)}$} to be the sheafification of this prestack in $\Perf/\Div^1$. 
\xqed{\blacktriangle}\end{dfn}

\begin{pro}\label{eqbetweentwisteddisplay} 
    Assume that for $i\not\in \Gamma$ we have $\mu_i=1$.
    Over $\Spd \mathbb{Z}_{p^r}$, there is an isomorphism $$\Xi_{(\CG,\boldsymbol{\mu},\Gamma)}\colon\CB^{\textup{perf}}_{\left(\CG,\boldsymbol\mu\right)}\to\CB^{\textup{perf}}_{\left(\CG^\Gamma,\mu^\Gamma\right)}$$ 
\end{pro}

\begin{proof}
    We know that $$\CB^{\textup{perf}}_{\left(\CG,\boldsymbol \mu\right)}=\CB^{\textup{perf,r}}_{\left(\CG_0,\mu\right)/\Div^r}$$ and we give an isomorphism from this stack to $\CB^{\textup{perf}}_{\left(\CG^\Gamma,\mu^\Gamma\right)}$. As usual, it is enough to define it on banal objects and morphisms. Define a map $\CG_0^r\to \CG_0^\Gamma$ by sending $\left(U_i\right)_{1\le i\le r}$ to $(\tilde{U}_{\gamma_t})$ where $$\tilde{U}_{\gamma_t}=U_{\gamma_{t+1}-1}\sigma\left(U_{\gamma_{t+1}-2}\right)\dots \sigma^{r_{t-1}-1}(U_{\gamma_{t}})$$ and define a map $\CG_0^{\boldsymbol{\mu}}\to \CG_0^{\Gamma,\mu}$ as the projection. An easy computation shows that these maps induce an isomorphism of quotient groupoids. 
\end{proof}

We use $\CB^{\textup{perf}}_{\left(\CG^\Gamma,\mu^\Gamma\right)}$ instead of $\CB^{\textup{perf}}_{\left(\CG,\mu\right)}$ because it gives us some ``empty'' space to perform our manipulations later on.

In what follows, it is a good idea to keep the following example in mind:

\begin{example}\label{mainexample}
    Let $r=3$, $\Gamma=\{1,3\}$, $\CG=\GL_3$ , $\mu_1=\left(1,z,z^2\right)$ and $\mu_3=\left(1,z,z\right)$. Then $\CG^{\mu_3}\left(A,d\right)$ consists of  matrices $\left(k_{ij}\right)$ where $\sigma\left(k_{13}\right)$ and $\sigma\left(k_{12}\right)$ are divisible by $d$, and $\CG^{\mu_1}\left(A,d\right)$ consists of  matrices $\left(k_{ij}\right)$ where $\sigma\left(k_{13}\right)$ and $\sigma\left(k_{23}\right)$ are divisible by $d$, while $\sigma\left(k_{13}\right)$ is divisible by $d^2$. In this setting, the action is given by the formula: $$\left(k_1,k_3\right).\left(b_1,b_3\right)=\left(k_1^{-1}b_1\Phi_d^{\mu_3}(k_3)\left(d\right)^{-1},k_3^{-1}b_3 \mu_1\Phi_{\sigma(d)}^{\mu_1}(\sigma(k_1))\right)$$
    The functor to  shtukas is given on objects by $$\textup{Fib}_{\GL_3}\left(\left(b_1,b_3\right)\right)=\left(b_1\mu_3(d),\sigma^{-1}\left(b_3\right)\mu_1(d),1\right)$$ and sends a morphism  $\left(k_1,k_3\right)$ to $\left(k_1,\sigma^{-1}\left(k_3\right),k_3\right)$.
\xqed{\blacksquare}\end{example}

\section{Definition of \texorpdfstring{$\Flex$}{Flex}}\label{definieflex}

B\"ultel introduced a functor from a class of graded higher displays, bounded by a non-minuscule cocharacter, to the category of displays. Using this functor, he constructed the special fiber of meta-unitary Shimura varieties, a class of Shimura varieties that contains some non-Abelian type Shimura varieties \cite{bultel2008pel}. The definition of $\Flex$ does not immediately generalize to mixed characteristics because the construction uses the fact that $p$ is Frobenius invariant in an essential way. In mixed characteristics, we have to replace $p$ with $V\left(1\right)$, which does not have this property. The construction of this functor has two main steps; in the first step, we give a functor from the graded higher displays, bounded by a non-minuscule cocharacter, to the shtukas with multiple legs bounded by a minuscule cocharacter. This functor is defined in a group-theoretic way, and we denote it by $\Flex^\triangleright$. In the second step, we combine the legs to obtain a shtuka with a single leg, bounded by a minuscule cocharacter.

In this section, we want to lift the map $\Flex$ to mixed characteristics. It helps us to understand the integral model of meta-unitary Shimura varieties. First, we have the following observation:

\begin{lem}\label{zarbi}
    Let $\mu$ and $\eta$ be two cocharacters of $\GL_n$, and $(A,d)$ an oriented prism. Then we have the following:
    \begin{enumerate}[label=(\roman*)] 
         \item $\GL_n^{\mu+\eta}(A,d)\subset \GL_n^{\mu}(A,d)\cap \GL_n^{\eta}(A,d)$ 
         \item $\Phi^{\mu}_d\left(\GL_n^{\mu+\eta}(A,d)\right)\subset \GL_n^{\eta}\left(A,\sigma(d)\right)$ 
         \item $\Phi^{\eta}_{\sigma\left(d\right)}\circ\Phi^{\mu}_d=\sigma\circ\Phi^{\mu+\eta}_{d}$ 
    \end{enumerate}
\end{lem}

Let us use this lemma to first define the map $\Flex^\triangleright$ in our toy  example:

\begin{example}
    Consider the setting  of Example \ref{mainexample}. Define the cocharacters $$\mu^\triangleright_1=\left(1,z,z\right),\, \mu^\triangleright_2=\left(1,1,z\right),\, \mu^\triangleright_3=\mu_3$$ Note that $\mu^\triangleright_1+\mu^\triangleright_2=\mu_1$. We want to define a map $$\CB^{\textup{perf}}_{\left(\GL_3^\Gamma,\mu^\Gamma\right)}\left(A,d\right)\to \CB^{\textup{perf, r}}_{\left(\GL_3^3,\mu^\triangleright\right)}\left(A,\left(d,d,\sigma\left(d\right)\right)\right)$$ It is enough to do this on banal objects and morphisms. We send a banal object $\left(b_1,b_3\right)$ to $\left(b_1,1,b_3\right)$ and a morphism $\left(k_1,k_3\right)$ between banal objects to $\left(k_1,\mu(d)\sigma(k_1)\mu(d)^{-1},k_3\right)$. Using Lemma \ref{zarbi}, one can easily see that this defines a map on quotient groupoids. It is easy to see that the associated shtukas are isomorphic to each other when we restrict them on $\mathcal{Y}_{[r,\infty]}$ for $r>\kappa\left(d\right)$.
\xqed{\blacksquare}\end{example}

\begin{dfn}\label{sparsedfn}
    We say that a display datum $\left(\GL_n^\Gamma,\mu^\Gamma\right)$ is \emph{sparse} if for all $i$, the highest weight of $\mu_{\gamma_i}$  is not larger than $r_i^+$.
\xqed{\blacktriangle}\end{dfn}

Now we want to generalize the construction in the above example. We use the notations defined in \ref{notationgamma}. 

\begin{construction}\label{tildemu}
    Let $\left(\GL_n^\Gamma,\mu^\Gamma\right)$ be a display datum and set ${\CG}^\triangleright\coloneqq\Res_{\BZ_{p^r}/\BZ_p} \GL_n$. From $\mu^\Gamma$ we obtain a cocharacter $\mu$ for ${\CG}^\triangleright$ defined over $\BZ_{p^r}$ where we set $\mu_i=1$ for $i\not\in \Gamma$. We are going to construct a cocharacter ${\boldsymbol{\mu}}^\triangleright=(\mu^\triangleright_j)_{j=1}^r$ of $\GL_n^r$ such that all $\mu^\triangleright_j$ are minuscule. Consider $n$-dimensional vector spaces $(V_i)_{i=1}^r$. For $\gamma\in\Gamma$, we have an action of $\BG_m$ on $V_{\gamma}$ via $\mu_\gamma$.  The action of $\BG_m$ on $V_{\gamma_i}$ by $\mu_{\gamma_i}$ gives us a grading $V_{\gamma_i}=\bigoplus_{h\in\BZ} V_{\gamma_i}^h$ where $\BG_m$ acts on $V_{\gamma_i}^h$ by weight $h$. For $\gamma_i\le j< \gamma_{i+1}$ and $h\in\BZ$, set $V_j^h\coloneqq V_{\gamma_i}^h$ and consider a new action where for $h<j-\gamma_i+1$, $\BG_m$ acts on $V_j^h$ by weight $1$ and for $h\ge j-\gamma_i+1$ it acts by weight $0$. We call the associated cocharacter $\mu^\triangleright_j$. One can check that for all $1\le j\le r$, $\mu^\triangleright_j$ is minuscule. A simple double counting argument shows that when $(\GL_n^\Gamma,\mu^\Gamma)$ is sparse, we have $\mu_{\gamma_i}=\sum_{\gamma_{i}\le j<\gamma_{i+1}}\mu^\triangleright_j$.
\xqed{\blacktriangledown}\end{construction}

\begin{rmk}
    Set ${\CG}^\triangleright\coloneqq\Res_{\BZ_{p^r}/\BZ_p} \GL_n$ and let $\mu$ and $\boldsymbol{\mu}^\triangleright$ be as in the Construction \ref{tildemu}. Let $B(G^\triangleright,{\mu})$ and $B(G,\boldsymbol\mu^\triangleright)$ be the Kottwitz sets (See \cite[Definition 7.2]{hedayatzadeh2024deformations} for the definition).
     Note that if $\mu$ is defined over $\BZ_{p^r}$ then we have
    \[
    \bar{\mu}\coloneqq\sum_{\tau \in \textup{Gal}(\BQ_{p^r}/\mathbb{Q}_p)} \mu^\tau = \sum_{\tau \in \textup{Gal}(\BQ_{p^r}/\mathbb{Q}_p)} \mu^{\triangleright,\tau}\eqqcolon\bar{\mu}^\triangleright
    \] Therefore, the sets $B({\CG}^\triangleright, \mu)$ and $B({\CG}^\triangleright, \boldsymbol\mu^\triangleright)$ are equal, as they only depend on the $\bar{\mu}$ and $\overline{\mu^\triangleright}$. By \cite[Proposition 7.4]{hedayatzadeh2024deformations} we know that these sets are the topological image of the maps $\CB^{\textup{perf}}_{{(\CG^\triangleright}, \mu)}\to \Bun_{\CG^\triangleright}$ and $\CB^{\textup{perf}}_{({\CG}^\triangleright, \boldsymbol\mu^\triangleright)}\to \Bun_{\CG^\triangleright}$, respectively. This provides further evidence that the stacks $\CB^{\textup{perf}}_{({\CG}^\triangleright, \mu)}$ and $\CB^{\textup{perf}}_{(\CG^\triangleright, \boldsymbol\mu^\triangleright)}$ are related.
\xqed{\lozenge}\end{rmk}

\begin{construction}
    Consider an sparse datum $\left(\GL_n^\Gamma,\mu^{\Gamma}\right)$ and the map $\Div^1\to \Div^r$ which sends $I$ to $\left(\sigma^{t^+\left(i\right)-1}\left(I\right)\right)_{i=1}^r$. We want to define $$\Flex^\triangleright\colon\CB^{\textup{perf}}_{\left(\GL_n^\Gamma,\mu^\Gamma\right)}\to \CB^{\perf, r}_{\left(\GL_n\vphantom{^\Gamma},\boldsymbol\mu^\triangleright\right)}$$ over this map. It is enough to define it for banal objects and morphisms between them. On banal objects, we define it via the natural injection $\GL_n^\Gamma\to \GL_n^r$. On morphisms between banal objects, we need to define it in a way that $$\Flex^\triangleright\left(\left(k_{\gamma_i}\right)\cdot\left(b_{\gamma_i}\right)\right)=\Flex^\triangleright\left(\left(k_{\gamma_i}\right)\right)\cdot\Flex^\triangleright\left(\left(b_{\gamma_i}\right)\right)$$ Using Lemma \ref{zarbi}  we should have $\Flex^\triangleright\left((k_{\gamma_i})\right)=(\tilde{k}_j)$ where $\tilde{k}_{\gamma_i}=k_{\gamma_i}$ and for $\gamma_i< j <\gamma_{i+1}$ we inductively define $$\tilde{k}_j=\Phi_{\mu^\triangleright_{j-1}}^{\sigma^{t^+\left(i-1\right)}(I)}\left(\tilde{k}_{j-1}\right)$$
\xqed{\blacktriangledown}\end{construction}

When we apply the functor $\Flex^\triangleright$, we introduce multiple legs. We prefer to have a functor to the category of shtukas with one leg, because shtukas with one leg have a more direct geometric interpretation: one can see them as the cohomology of motives. We have multiple legs for $\textup{Fib}_{\CG}\left(\Flex^\triangleright\left(\CD\right)\right)$ but all these legs are Frobenius translates of each other, suggesting one could ``combine'' these legs to form a shtuka with one leg. This is what we want to do next; to do this, we need to view shtukas as modifications of vector bundles. Let us recall this perspective:

\begin{dfn}\label{dataofmodification}
    Assume that $y$ is a Cartier divisor on $\mathcal{Y}_S$ and $I_y$ to be the corresponding ideal. We denote the completion of $\mathcal{O}_{\mathcal{Y}_S}$ by $B^+_{\Div_{\mathcal{Y}}}(y)$ and $B^+_{\Div_{\mathcal{Y}}}(y)[\frac{1}{I_y}]$ by $B_{\Div_{\mathcal{Y}}}(y)$. The \emph{data of modification} of a $\CG$-torsor $\SP$ on $\CY_S$ at $y$ is a $G$-torsor $\SP'_y$ over $B^+_{\Div_{\mathcal{Y}}}(y)$ and an isomorphism of $\CG$-torsors $\SP'_y\left[\frac{1}{I_y}\right]\xrightarrow{\sim}\SP_y\left[\frac{1}{I_y}\right]$ over $B^+_{\Div_{\mathcal{Y}}}(y)$. 
\xqed{\blacktriangle}\end{dfn}

\begin{rmk}\label{Bouvillelaslo}
    By Beauville--Laszlo theorem \cite{zbMATH00810980}, having a data of modification for $\SP$ at $y$ is equivalent to having a $\CG$-torsor $\SP'$ and a morphism $\SP\to \SP'$ that is an isomorphism away from $y$. 
\xqed{\lozenge}\end{rmk}

There is an equivalence between the category of $\CG$-shtukas over $\mathcal{Y}_S$ with no legs and the category of $\CG(\mathbb{Z}_p)$-torsors on $S$. 
Consider the multi-leg $\CG$-shtuka $\left(\SP,y_1,...,y_n,\eta\right)$ over $S$. Assume that $\epsilon<\min\{\kappa(y_i)\}_{i=1}^n$. Consider the restriction of $\SP$ to $\mathcal{Y}_{S,\left(0,\epsilon\right)}$. This is a $\CG$-shtuka over $\CY_{S,(0,\epsilon)}$ with no leg, and because Frobenius is a contraction (with respect to $\kappa$), one can extend it to a $\CG$-shtuka over the whole curve $\CY_S$. We denote this shtuka by $\SP_0$. Therefore, we obtain a $G\left(\mathbb{Z}_p\right)$-torsor $T$ on $S$. By \cite[Section 12]{scholze2020berkeley}, having the shtuka $\SP$ is equivalent to having a shtuka with no leg$\SP_0$ and the data of modifications of $\SP_0$ at $y_1,y_2,\dots,y_n$. One can recover $\SP$ from $\SP_0$ by first modifying $\SP_0$ at $y_0$, then modifying the resulting shtuka at $y_1$, and so on. 

One can ``combine'' legs to obtain a $G$-shtuka with one leg. Let us first set up a notation.

\begin{dfn}
    Consider the locus inside $\Div^r$ where the untilts are given by the ideals $y,\sigma(y),\dots,\sigma^r(y)$, we denote the restriction of $\Sht_{(\CG,\mu)}$ to this locus by $\Sht_{\left(\CG,\mu,\left(y,\sigma\left(y\right),\cdots,\sigma^r\left(y\right)\right)\right)}$.
\xqed{\blacktriangle}\end{dfn}

\begin{pro}\label{com}
    There is a canonical functor $$\com\colon\Sht_{\left(\CG,(\mu_i),(y,\sigma(y),\cdots,\sigma^r(y))\right)}\to \Sht_{(\CG,\sum \mu_i,y)}$$
\end{pro}

\begin{proof}
    let $\SP\in \Sht_{\left(\CG,\left(\mu_i\right),\left(y,\sigma\left(y\right),\cdots,\sigma^r\left(y\right)\right)\right)}$ be a $\CG$-shtuka. By the above paragraph, it provides us with a $\CG(\mathbb{Z}_p)$-torsor $T$ and the data of modifications at $y,\sigma\left(y\right),\dots,\sigma^r\left(y\right)$. Let $\SP_0$ be the $\CG$-shtuka with no leg associated with $\SP$. The shtukas $\SP_0$ and $\SP$ are isomorphic to each other in a punctured neighborhood of $y$. We define $\com(\SP)$ as the shtuka with one leg, obtained from $\SP_0$ and the data of modification coming from the isomorphism between $\SP$ and $\SP_0$ near $ y$.
\end{proof}

\begin{rmk}\label{goingbetweendisplaysandshtukas}
    Note that $\SP$ and $\com(\SP)$ are isomorphic to each other near infinity, therefore, if one can extend $\SP$ to infinity one can also extend $\com (\SP)$ to infinity. By \cite[Corollary 5.6.10]{ito2023prismatic} we have an equivalence of category between the $(\CG,\mu)$-displays over $(R,R^+)$ and the $(\CG,\mu)$-displays over $R^+$. By \cite[Proposition 9.2]{hedayatzadeh2024deformations} there is an equivalence of category between $(\CG,\mu)$-displays over $R^+$ and Shtukas over $R^+$. Therefore, if the Shtuka $\SP$ over $R^+$ comes from a $(\CG,\mu)$-display $\CD$, then $\com (\SP)$ comes from a $(\CG,\mu)$-display $\com(\CD)$. When one perform the modification using Beauville--Laszlo theorem, the cocharacters sum up, therefore we obtain a shtuka bounded by $\sum \mu_i$.
\xqed{\lozenge}\end{rmk}

\begin{rmk}
    By definition, the shtuka $\com\left(\SP\right)$ is isomorphic to $\SP$ on $\mathcal{Y}_{S,(\kappa\left(y\right),\infty)}$. By \cite{bartling2022mathcal}, this isomorphism gives an isogeny between the reduction of the associated displays modulo $p$.
\xqed{\lozenge}\end{rmk}

\begin{dfn}\label{defFLex}
    Set $\CG^\triangleright\coloneqq \Res_{\mathbb{Z}_{p^r}/\mathbb{Z}_p}\GL_n$. We define the functor $$\Flex\colon\CB^{\textup{perf}}_{(\CG^\triangleright,\mu)}\to \Sht^{\textup{perf}}_{(\CG^\triangleright,\mu^\triangleright)}$$ to be the composition $\com\circ\textup{Fib}\circ \Flex^\triangleright\circ\Xi_{(\CG,\boldsymbol{\mu},\Gamma)}$. Using remark \ref{goingbetweendisplaysandshtukas}, we upgrade this functor to the functor $$\Flex\colon\CB^{\textup{perf}}_{(\CG^\triangleright,\mu)}\to \CB^{\textup{perf}}_{(\CG^\triangleright,\mu^\triangleright)}$$
\xqed{\blacktriangle}\end{dfn}

\begin{rmk}\label{isogenyrighttriangle}
    There is an isogeny away from legs between $\CD$ and $\Flex^\triangleright\left(\CD\right)$. By descent, we only need to define it for banal displays; for them, we define the isogeny by $\left(U_i\right)_{i=0}^r$ where $U_i=1$ for $i\in \Gamma$ and $$U_i=\prod_{\gamma_t<j\le i}\mu^\triangleright_{j-1}\left(\sigma^{t^+\left(i\right)-1}\left(d\right)\right)$$
for $\gamma_t<i<\gamma_{t+1}$.
\xqed{\lozenge}\end{rmk}

\begin{dfn}
    An \emph{isogeny modulo $p$} between shtukas $\SP$ and $\SP'$ is an isomorphism $i_t\colon\SP_{[t,\infty)}\to \SP_{[t,\infty)}$ for some $t>0$.  By \cite[Lemma 8.1]{bartling2022mathcal}, it is equivalent to having an isogeny between the reduction mod $p$ of $\SP$ and $\SP'$.
\xqed{\blacktriangle}\end{dfn}

\begin{pro}\label{isogenyclass}
    Let $S$ be a perfectoid space and $\CD\in \CB^{\textup{perf}}_{(\CG_0^\Gamma,\mu^\Gamma)}(S)$. There is a canonical isomorphism between $\textup{Fib}\left(\CD\right)_{\mathcal{Y}_{S,\left[t,\infty\right)}}$ and $\Flex\left(\CD\right)_{\mathcal{Y}_{S,\left[t,\infty\right)}}$ for big enough $t\in \BR_{>0}$. In particular, we obtain an isogeny modulo $p$ between the associated displays.
\end{pro}

\begin{proof}
    We saw before that there is a canonical isogeny between $\CD$ and $\Flex^\triangleright\left(\CD\right)$, and there is an isomorphism between $\textup{Fib}(\CD)_{\left(\kappa\left(y\right),\infty\right)}$ and $\com\left(\CD\right)_{\left(\kappa\left(y\right),\infty\right)}$. 
\end{proof}

\begin{dfn}
    Let $\CG$ be a reductive group over $\mathbb{Z}_p$, $\mu$ a cocharacter for $\CG$ and $b\in B(\CG)=\CG(\breve{\BQ}_p)/\sigma\textup{-conj}$ such that $(\CG,\mu,b)$ is an integral Shimura datum in the sense of \cite{bueltel2020g} (see also Appendix \ref{shtukaappaendix}). We define the banal display $\bar{\CD}_b$ as the banal display represented by any element of $\CG(\breve{\BQ}_p)$ lifting $b\mu(p)^{-1}$.
\xqed{\blacktriangle}\end{dfn}

\begin{dfn}
    Let $(\CG,\mu,b)$ be an integral Shimura datum. We define the stack $\CB^{\textup{perf}}_{(\CG,\mu,b)}$ as follows: for an affinoid perfectoid space $\Spa (R,R^+)$, objects of $\CB^{\textup{perf}}_{(\CG,\mu,b)}\left(\Spa (R,R^+)\right)$ are pairs $(\CD,\eta)$ with $\CD\in \CB^{\textup{perf}}_{({\CG},\mu)}(R,R^+)$ and $\eta$ is a quasi-isogeny $\eta\colon\bar{\CD}\dashrightarrow \bar{\CD}_b$ where $\bar{\CD}$ is the reduction of $\CD$ mod $p$. A morphism from $(\CD,\eta)$ to $(\CD',\eta')$ is a morphism $k\colon\CD\to \CD'$ in $\mathcal{B}^{\textup{perf}}_{(\CG,\mu)}$ such that $\eta'\circ k=\eta$. 
\xqed{\blacktriangle}\end{dfn}

\begin{construction}
    We want to lift $\Flex$ to a functor $$\Flex\colon\CB^{\textup{perf}}_{(\CG,\mu,b)}\to \CB^{\textup{perf}}_{({\CG^\triangleright},\mu^\triangleright,b)}$$
    Let $\left(\CD,\eta\right)\in\CB^{\textup{perf}}_{(\CG,\mu,b)}$. Define $\Flex\left(\CD,\eta\right)$ to be $\left(\Flex\left(\CD\right),\Flex\left(\eta\right)\right)$ where $\Flex\left(\eta\right)$ is the isogeny obtained from composing $\eta$ with the isogeny between $(\bar{\CD})$ and $\overline{\Flex(\CD)}$ given in Proposition \ref{isogenyclass}.
\xqed{\blacktriangledown}\end{construction}

Over the generic fiber, one can define level structures on the stack of shtukas. Let $K\subset \CG(\mathbb{Z}_p)$ be an open compact subgroup. By a $K$-level structure on a $\CG$-shtuka $\SP$, we mean a $K$-torsor inside the associated $\CG(\mathbb{Z}_p)$-torsor $T$. Since $\SP$ and $\Flex\left(\SP\right)$ are isomorphic near $\mathcal{Y}_{S,[0,0]}$ (from now on, when there is no ambiguity, we denote $\mathcal{Y}_{S,[0,0]}$ by $0$), one can extend the definition of $\Flex$ to an arbitrary level:

\begin{dfn}\label{Flexarbirarylevel}
    Let $K_p\subset \GL_n(\mathbb{Z}_p)$ be a compact open subgroup. Define the map 
    $$\Flex\colon\Sht_{(G,\mu,b),K_p,\mathbb{Q}_p}\to \Sht_{(G^\triangleright,\mu^\triangleright,b),K_p,\mathbb{Q}_p}$$ by sending $\left(\SP,\tau,T\right)$ to $\left(\textup{Fib}( \Flex\left(\SP\right)),\Flex\left(\tau\right),T\right)$.
\xqed{\blacktriangle}\end{dfn}

Let us take a closer look at the functor $\Flex$ defined above. Usually, we have a better bound than $\sum \mu^\triangleright_i$ at the poles, and in some examples, we obtain a shtuka bounded by a minuscule cocharacter.

\begin{example}
    Consider the datum $\Gamma=\{1\}\subset \{1,2,3\}$ and $\mu_1=\left(1,z^3\right)$, and let $y$ be a Cartier divisor on the curve $\mathcal{Y}_S$. Let $\CD\in \CB^{\textup{perf}}_{\left(\CG^\Gamma,\mu^\Gamma\right)}$ and $\SP=\Fib(\Flex\left(\CD\right))$. If $$\SP=\SP_1\oplus \SP_2\oplus \SP_3$$, we have a morphism $\sigma^* \SP_i\to \SP_{i+1}$ that is an isomorphism away from $\sigma^{i-1}\left(y\right)$ and is bounded by the cocharacter $\mu\coloneqq \left(1,z\right)$ at this point. To compute $\com\left(\SP\right)$, we consider the shtuka with no legs $\SP^0$ that extends $\SP_{\left(0,\kappa\left(\sigma^3\left(y\right)\right)\right)}$ and compare $\SP_y$ and $\SP^0_y$. We have an isomorphism $\SP^0_{\sigma^3\left(y\right)}\to \SP^0_y$ and the morphisms $$\SP_{i,\sigma^3\left(y\right)}\to \SP_{i-2,\sigma^2\left(y\right)}\to \SP_{i-1,\sigma\left(y\right)}\to \SP_{i,y}$$ For $i=2$, only the third morphism is not an isomorphism; for $i=1$, only the second one is not an isomorphism; and for $i=3$, only the third one is not an isomorphism. These morphisms are all bounded by $\mu$. 

Therefore, when we combine the legs we obtain a shtuka that is bounded by $\mu$ at each grade. For convenience, let us summarize this information in a table.

\vspace{0.5 cm}
\begin{center}
\resizebox{\textwidth}{!}{
\begin{tikzpicture}
\node at (0,0) {
\begin{tabular}{|c|l|}
\hline
$\sigma^*\SP_1 \to \SP_2$ & isomorphism away from $y$, modification data $\Xi_1$ at $y$\\
\hline
$\sigma^*\SP_2 \to \SP_3$ & isomorphism away from $\sigma(y)$, modification data $\Xi_2$ at $\sigma(y)$\\
\hline
$\sigma^*\SP_3 \to \SP_1$ & isomorphism away from $\sigma^2(y)$, modification data $\Xi_3$ at $\sigma^2(y)$\\
\hline
\end{tabular}
};
\end{tikzpicture}
}
\end{center}

\begin{center}
\resizebox{\textwidth}{!}{
\begin{tikzpicture}
    \coordinate (table_top_left) at (0,0);

    \pgfmathsetmacro{\mycolAwidth}{2}  
    \pgfmathsetmacro{\mycolBwidth}{2.7}  
    \pgfmathsetmacro{\mycolCwidth}{2.7}  
    \pgfmathsetmacro{\mycolDwidth}{3.9}  
    \pgfmathsetmacro{\myrowheightone}{2.0} 
    \pgfmathsetmacro{\myrowheighttwo}{1.5} 

    \pgfmathsetmacro{\mytablewidth}{\mycolAwidth + \mycolBwidth + \mycolCwidth + \mycolDwidth}
    \pgfmathsetmacro{\mytableheight}{\myrowheightone + \myrowheighttwo}

    \draw (table_top_left) rectangle ++(\mytablewidth, -\mytableheight);

    \draw (table_top_left) ++(0, -\myrowheightone) -- ++(\mytablewidth, 0);

    \draw (table_top_left) ++(\mycolAwidth, 0) -- ++(0, -\mytableheight);
    \draw (table_top_left) ++(\mycolAwidth + \mycolBwidth, 0) -- ++(0, -\mytableheight);
    \draw (table_top_left) ++(\mycolAwidth + \mycolBwidth + \mycolCwidth, 0) -- ++(0, -\mytableheight);

    \node at ($(table_top_left) + (0.5*\mycolAwidth, -0.5*\myrowheightone)$) {$\SP_{1,y},\SP_{3,\sigma(y)}$};
    \node at ($(table_top_left) + (\mycolAwidth + 0.5*\mycolBwidth, -0.5*\myrowheightone)$) {$\SP_{3,\sigma(y)},\SP_{2,\sigma^2(y)}$};
    \node at ($(table_top_left) + (\mycolAwidth + \mycolBwidth + 0.5*\mycolCwidth, -0.5*\myrowheightone)$) {$\SP_{3,\sigma^2(y)},\SP_{1,\sigma^3(y)}$};
    \node at ($(table_top_left) + (\mycolAwidth + \mycolBwidth + \mycolCwidth + 0.5*\mycolDwidth, -0.5*\myrowheightone)$) {$\tilde{\SP}_{1,d} ,\tilde{\SP}_{1,\sigma^3(y)},\SP_{1,\sigma^3(y)}$};

    \node[align=center] at ($(table_top_left) + (0.5*\mycolAwidth, -\myrowheightone - 0.5*\myrowheighttwo)$) {iso};
    \node[align=center] at ($(table_top_left) + (\mycolAwidth + 0.5*\mycolBwidth, -\myrowheightone - 0.5*\myrowheighttwo)$) {mod. data\\$\Xi_2$};
    \node at ($(table_top_left) + (\mycolAwidth + \mycolBwidth + 0.5*\mycolCwidth, -\myrowheightone - 0.5*\myrowheighttwo)$) {iso};
    \node at ($(table_top_left) + (\mycolAwidth + \mycolBwidth + \mycolCwidth + 0.5*\mycolDwidth, -\myrowheightone - 0.5*\myrowheighttwo)$) {iso};

    \pgfmathsetmacro{\myarrowxstart}{\mytablewidth + 0.5} 
    \pgfmathsetmacro{\myarrowxend}{\myarrowxstart + 0.5} 
    \pgfmathsetmacro{\myarrowycenter}{-\mytableheight/2} 

    \draw[->, line width=0.8pt] (\myarrowxstart, \myarrowycenter) -- (\myarrowxend, \myarrowycenter);

    \pgfmathsetmacro{\myboxwidth}{2.0} 
    \pgfmathsetmacro{\myboxheight}{3.5} 
    
    \pgfmathsetmacro{\myboxxstart}{\myarrowxend + 0.5}
    \pgfmathsetmacro{\myboxycenter}{-\mytableheight/2}
    \pgfmathsetmacro{\myboxytopedge}{\myboxycenter + 0.5*\myboxheight}
    \pgfmathsetmacro{\myboxxcenter}{\myboxxstart + 0.5*\myboxwidth}

    \draw (\myboxxstart, \myboxytopedge) rectangle ++(\myboxwidth, -\myboxheight);
    \pgfmathsetmacro{\myboxseparatorY}{\myboxytopedge - 0.5 * \myboxheight}
    \draw (\myboxxstart, \myboxseparatorY) -- (\myboxxstart + \myboxwidth, \myboxseparatorY);

    \node[align=center] at (\myboxxcenter, {(\myboxytopedge + \myboxseparatorY) / 2}) {
        $\tilde{\SP}_{1,y} \SP_{1,y}$
    };

    \node[align=center] at (\myboxxcenter, {(\myboxseparatorY + (\myboxytopedge - \myboxheight)) / 2}) {
        mod. data\\$\Xi_2$
    };
\end{tikzpicture}
}
\end{center}
\xqed{\blacksquare}\end{example}

More generally, assume that $m$ is an odd integer, and for each $i$, we have $r^+_i=m$. Similar to the above, at each grade, only one modification occurs, and when we combine the legs, we obtain a modification that is bounded by a minuscule cocharacter at each grade.

\begin{dfn}\label{regular}
    We call a sparse datum $(\GL_n^\Gamma,\mu^\Gamma)$ \emph{regular} if for all $w\in \BZ/r$, there is at most one $w'\in\BZ/r$ such that $\textup{mod}(w-w')=t^+ (w'-1)$. 
\xqed{\blacktriangle}\end{dfn}

\begin{rmk}
\begin{enumerate}[label=\arabic*)]
    \item Assume that $(\GL_n^\Gamma,\mu^\Gamma)$ is sparse. If all intervals $[\gamma_i,\gamma_{i+1})$ (for all $\gamma_i\in\Gamma)$ have the same length and that length is odd then the datum is regular. However, there are cases where the datum is regular, but these intervals are not all equal of odd length.
    \item Now, assume that the datum is regular, then the underlying shtuka after applying $\Flex$ (so, after combining the legs) is bounded by a minuscule cocharacter. 
    \item Assume again that $(\GL_n^\Gamma,\mu^\Gamma)$ is regular. Then $\com$ is a monomorphism, because when we combine the legs, modifications occur at separate grades and therefore we can recover the modifications needed to obtain $\SP$ from the ones for $\com(\SP)$.
\end{enumerate}
\xqed{\lozenge}\end{rmk}
   
\begin{example}
    Some sparse data are not regular. For example, if $$\Gamma=\{1,4\}\subset \{1,2,3,4\}$$ and $\mu_1=\left(1,z^4\right),\mu_4=\left(1,z\right)$, then $\Flex(\CD) $ is bounded by $\left(1,z^2\right)$ at grades $1,3$, by $\left(1,z\right)$ at grade $4$, and is an isomorphism at grade $2$.
\xqed{\blacksquare}\end{example}

\section{Sparse Local Shimura datum}\label{sectionsparse}

Let $\CG_0$ be a reductive group over $\mathbb{Z}_{p^r}$. Let $\mu$ be a cocharacter of $\CG\coloneqq {\Res}_{\mathbb{Z}_{p^r}/\mathbb{Z}_p} \CG_0$ defined over $\BZ_q$. We have an isomorphism $\CG_{\mathbb{Z}_{p^r}}\xrightarrow{\sim}
\CG_0^r$. From a representation $\left(V,\rho\right)$ of $\CG_0$ over $\mathbb{Z}_{p^r}$, we obtain a map $\CG_{\mathbb{Z}_{p^r}}\to \GL\left(V\right)^{r}$. By composing this map with $\mu$, we obtain a cocharacter $\left(\mu_{\tau}\right)_{i=1}^{r}$ of $\GL\left(V\right)^r$. Let $\Gamma=\{\gamma\in \mathbb{Z}/r\mathbb{Z}\,|\,\mu_\gamma\not=1\}$ and consider the display datum $\left(\GL\left(V\right)^\Gamma,\mu^\Gamma\right)$.

\begin{dfn}\label{sparseShimuradfn}
    In the above setting, we say that an integral Shimura datum $(\CG,\mu,b)$ is \emph{sparse} if $\CG_0$ has a faithful representation $(V,\rho)$ such that $\left(\GL\left(V\right)^\Gamma,\mu^\Gamma\right)$ is sparse. We say that it is \emph{regularly sparse} if $\left(\GL\left(V\right)^\Gamma,\mu^\Gamma\right)$ is also regular in the sense of Definition \ref{regular}. 
\xqed{\blacktriangle}\end{dfn}

Let $\CG_0\subset \GL(V)$ be the stabilizer of tensors $\left(s_\alpha\right)$ of $V$. Set ${\CG}^\triangleright_\rho\coloneqq \textup{Res}_{\mathbb{Z}_{p^r}/\mathbb{Z}_p} \GL\left(V\right)$. Let us still denote by $s_\alpha$ the tensor of $\oplus_{i=1}^r V$ given by the diagonal embedding of $V\into\oplus_{i=1}^r V$. Then, $\CG_{\mathbb{Z}_{p^r}}$ is the stabilizer of tensors $(s_{\alpha})$  inside ${\CG}^\triangleright_{\rho, \mathbb{Z}_{p^r}}=\GL(V)^r$. By a slight abuse of notation, we denote by $\Flex$ the composition of $\Flex$ with the natural morphism $\CB^{\textup{perf}}_{(\CG,\mu)}\to\CB^{\textup{perf}}_{(\CG_\rho^\triangleright,\mu)}.$

The functor \(\Flex\colon \CB^{\perf}_{(\CG, \mu)} \to \CB_{(\CG_\rho^\triangleright, \mu^\triangleright)}^{\textup{perf}}\) is faithful but not full. To achieve fullness, we need to restrict the morphisms on the target category to those that respect the \(\CG\)-structure in some sense, for instance, those that preserve the tensors \(s_{\alpha}\). However, \(\Flex\) is not an isogeny; it is only an isogeny modulo \(p\). Consequently, we cannot directly transfer the tensors \((s_{\alpha})\) from \(\CD\) to \(\Flex(\CD)\), but we can do so modulo \(p\). Therefore, we restrict to morphisms that preserve these tensors modulo \(p\). This restriction suffices to ensure the fullness of \(\Flex\), since taking displays modulo \(p\) is a faithful functor. This motivates the following definition:

\begin{dfn}\label{tanakianversionofB}
    Let \emph{$\CB_{(\CG_\rho^\triangleright,\mu^\triangleright)}^{\textup{tan}}$} be the sheaf of groupoids over $\Spd\mathbb{Z}_q$, sending $\Spa\left(R,R^{+}\right)\to \Spd \BZ_q$, given by an untilt $\Spa(R^{\sharp},R^{\sharp+})$, to pairs $\left(\SD,\left(t_\alpha\right)\right)$ such that $\SD\in \CB^{\textup{perf}}_{(\CG_\rho^\triangleright,\mu^\triangleright)} $ and $\left(t_\alpha\right)$ are Frobenius-invariant tensors in $\left(\SD_{A_{\textup{crys}}\left(R^{\sharp+}\right/p})\right)^\otimes$ and \'etale locally there exists an isomorphism of graded vector bundles $$\SD\otimes A_{\textup{crys}}\left(R^{\sharp+}/p\right)\to V\otimes_{\mathbb{Z}_p}  A_{\textup{crys}}\left(R^{\sharp+}/p\right)$$ that sends $\left(t_\alpha\right)$ to $\left(s_\alpha\right)$.
    The morphisms are the morphisms between displays that preserve tensors.    
\xqed{\blacktriangle}\end{dfn}
 
\begin{rmk}\label{FlexRdfn}
    One can upgrade $\Flex\colon\CB^{\textup{perf}}_{(\CG,\mu)}\to \CB^{\textup{perf}}_{(\CG_\rho^\triangleright,\mu^\triangleright)}$ to a map $\Flex^{\textup{tan}}\colon\CB_{(\CG, \mu)}^{\textup{perf}}\to \CB_{(\CG_\rho^\triangleright,\mu^\triangleright)}^{\textup{tan}}$ as follows. By Proposition \ref{isogenyclass} there is an isogeny $\eta$ modulo $p$ between  $\CD$ and $\Flex(\CD)$. By Tannakian view point on $G$-torsors, we have the tensors $(s_\alpha)$ inside $\CD\otimes A_{\textup{crys}}\left(R^{\sharp,+}/p\right)$, we define the tensors $t_\alpha\coloneqq \eta(s_{\alpha})$ inside $\Flex(\CD)\otimes A_{\textup{crys}}\left(R^{\sharp,+}/p\right)$.
\xqed{\lozenge}\end{rmk}

\begin{pro}\label{mono}
    The map $\Flex^{\textup{tan}}\colon\CB_{(\CG, \mu)}^{\textup{perf}}\to \CB_{(\CG_\rho^\triangleright,\mu^\triangleright)}^{\textup{tan}}$ is fully faithful.
\end{pro}

\begin{proof}
    The map $\Flex^{\textup{tan}}$ is the composition of four maps: the first map is an equivalence of categories. The second map is induced from the closed immersion $\rho\colon \CG\to \GL\left(V\right)$ and the functoriality. It is clear that it is faithful. To see that it is full, assume we have maps $\tilde{\CD}=\rho\left(\CD\right)\to \tilde{\CD'}=\rho\left(\CD'\right)$ that stabilize tensors $\left(s_\alpha\right)$ on $\textup{Fib}\left(\tilde{\CD}\right)\otimes A_{\textup{crys}}\left(R^{\sharp,+}\right)$. After going to a suitable cover, we can assume these maps are given by $$k\in \GL\left(V\otimes A_{\textup{inf}}(R^+)\right)\cap \CG\left(A_{\textup{crys}}\left(R^{\sharp,+}/p\right)\right)$$ And we have $$\CG\left(A_{\textup{inf}}(R^{\sharp,+})\right)=\CG\left(A_{\textup{crys}}\left(R^+/p\right)\right)\cap \GL\left(V\otimes A_{\textup{inf}}(R^+)\right)$$ Therefore, $k$ comes from a morphism between $\CD$ and $\CD'$. 
    
    The third map is $\Flex^\triangleright$. We want to prove that this map is fully faithful. By descent, it is enough to prove this for banal objects and morphisms. From the definition, it is clear that this map is faithful, so it is enough to show that it is full. To prove this, note that there is an isomorphism between $\CD$ and $\Flex D$ away from $\{\sigma^i\left(d\right)\}$, given by a matrix $T$. Therefore, an isomorphism $k\colon\Flex^\triangleright\left(\CD\right)\xrightarrow{\sim} \Flex^\triangleright\left(\CD'\right)$ induces an isomorphism $$k\colon\CD\left[{\frac{1}{\prod\sigma^i\left(d\right)}}\right]\xrightarrow{\sim} \CD'\left[{\frac{1}{\prod \sigma^i\left(d\right)}}\right]$$ Consider the lattices $\CD$ and $\CD'$ inside $\CD\left[{\frac{1}{\prod \sigma^i\left(d\right)}}\right]$ and $\CD'\left[{\frac{1}{\prod \sigma^i\left(d\right)}}\right]$. We have

    $$k\left(\CD\right)=k\left(T\cdot\Flex\left(\CD\right)\right)=T\cdot k\left(\Flex\left(\CD\right)\right)=T\cdot\Flex\left(\CD'\right)=\CD'$$

    Thus, $k$ comes from an isomorphism between $\CD$ and $\CD'$. The final map is $\com$, which is also fully faithful by Proposition \ref{com}.
\end{proof}

\section{Regularly Sparse Local Shimura Varieties}\label{regularlysparsesection}

\begin{dfn}
    We call the local Shimura varieties associated with regularly sparse Shimura data \emph{regularly sparse local Shimura varieties}. Let $(\CG,\mu,b)$  be an integral regularly sparse Shimura datum (cf. Definition \ref{sparseShimuradfn}) and $\mathcal{M}^{\textup{int}}_{(\CG,\mu,b)}$ be the associated integral local Shimura variety. We will denote the generic fiber of the integral local Shimura variety by $\CM_{(\CG,\mu,b)}$ For the definition and basic properties of local Shimura varieties, see Appendix \ref{shtukaappaendix} and \cite{pappas2022p}. By definition, for every perfectoid space $(S,S^+)$ we have a bijection between $\mathcal{M}^{\textup{int}}_{(\CG,\mu,b)}(S,S^+)$ and the set of isomorphism classes of objects in $\Sht_{(\CG,\mu,b)}(S,S^+)$, therefore we have the map $$\Flex\colon\mathcal{M}^{\textup{int}}_{(\CG,\mu,b)}\to \mathcal{M}^{\textup{int}}_{({\CG}^\triangleright_\rho,\mu^\triangleright,b)}$$
\xqed{\blacktriangle}\end{dfn}

 In this section, we want to prove Conjectures \hyperlink{representintegrallocal}{A} and \hyperlink{perfectoidintegrallocal}{B} for regularly sparse local Shimura varieties. First, we prove that the map $\Flex$ is a closed immersion between the integral models of local Shimura varieties.

\begin{rmk}
    Note that for a perfectoid ring \((R, R^+)\), there is a natural bijection
    \[
    \textup{Fib}_{\CG} \colon \mathcal{B}^{\textup{perf}}_{(\CG, \mu, b)}(R, R^+)\to \mathcal{M}^{\textup{int}}_{(\CG, \mu, b)}(R, R^+)
    \]
    (see \cite[Proposition 8.2]{hedayatzadeh2024deformations}).
    We will frequently use this correspondence in what follows to pass between displays and points of the local Shimura variety.
\xqed{\lozenge}\end{rmk}

\begin{lem}\label{Flexquasicomact}
    Let $\left(\CG,\mu,b\right)$ be a regularly sparse Shimura datum. The map $$\Flex\colon\mathcal{M}^{\textup{int}}_{(\CG,\mu,b)}\to \mathcal{M}^{\textup{int}}_{({\CG}^\triangleright_\rho,\mu^\triangleright,b)}$$ is quasi-compact.
\end{lem}

\begin{proof}
    The proof is similar to that of \cite[Proposition 3.6.2]{pappas2022p}. Let $(R,R^+)$ be an integral perfectoid ring and let $d$ be a generator of $\kernel \theta_{R^+}$. For $F\in \CG\left(W\left(R^+\right)[\frac{1}{d}]\right)$, let $\SP_{F}$ be the $\CG$-shtuka with a leg at $d$ with the trivial underlying vector bundle and Frobenius $F$.   Consider the v-sheaf $\mathcal{LM}^{\textup{int}}_{(\CG,\mu,b)}$ which maps $R^+$ to pairs $\left(F,i_r\right)$ where $F\in \CG\left(W\left(R^+\right)[\frac{1}{d}]\right)$ and $i_a$ is an isomorphism between $\SP_{F,[a,\infty]}$ and $\SP_{b,[a,\infty]}$ for some $a>0$, and we require that $\SP_{F}$ be bounded by $\mu$.
    
     This v-sheaf is a v-cover of $\mathcal{M}^{\textup{int}}_{(\CG,\mu,b)}$. Therefore, it is enough to prove that one can lift $\Flex^\triangleright$ to a quasi-compact map $$\mathcal{L}\Flex^\triangleright\colon\mathcal{LM}^{\textup{int}}_{(\CG,\mu,b)}\to \mathcal{LM}^{\textup{int}}_{({\CG}^\triangleright_\rho,\mu^\triangleright,b)}$$ Over this cover, $\mathcal{L}\Flex^\triangleright$ is the natural inclusion in the first component and the composition with the canonical quasi-isogeny between $F$ and $\Flex^\triangleright(F)$ in the second component. This map is quasi-compact, the proof being identical to that of \cite[Proposition 3.6.2]{pappas2022p}.
\end{proof}

\begin{lem}\label{colsedimmersionofvsheaf}
    Let $\left(\CG,\mu,b\right)$ be a regularly sparse Shimura datum. The map $$\Flex\colon\mathcal{M}^{\textup{int}}_{(\CG,\mu,b)}\to \mathcal{M}^{\textup{int}}_{({\CG}^\triangleright_\rho,\mu^\triangleright,b)}$$ is a closed immersion of v-sheaves.
\end{lem}

\begin{proof}
    This map is quasi-compact between partially proper spaces; hence, it is proper. By Proposition \ref{mono}, the map $\Flex^{\textup{tan}}$ is fully faithful, and thus induces an injection from the set of isomorphism classes of objects in $\CB_{(\CG,\mu,b)}$ to the set of isomorphism classes of objects in $\CB_{(\CG_\rho^\triangleright,\mu,b)}^{\textup{tan}}$. Moreover, there is an injection from the set of isomorphism classes of objects in $\CB_{(\CG,\mu,b)}^{\textup{tan}}$ to the set of isomorphism classes of objects in $\CB_{({\CG},\mu^\triangleright,b)}$. Therefore, the map \[\Flex\colon\mathcal{M}^{\textup{int}}_{(\CG,\mu,b)}\longrightarrow \mathcal{M}^{\textup{int}}_{({\CG}^\triangleright_\rho,\mu^\triangleright,b)}\] is a proper injection, which implies that it is a closed immersion.
\end{proof}

 Let us recall the definition of the specialization map for v-sheaves.
 
\begin{dfn}
    Let $\mathcal{F}$ be a small v-sheaf. We say that $\mathcal{F}$ is \emph{specializing} if the following conditions hold:
    \begin{enumerate}
        \item[\textup{(1)}] \textbf{v-Locally formal:}
        There exists a family $\{(B_i,B_i^+)\}_{i\in I}$ of formal Huber pairs over $\mathbb{Z}_p$ and a surjective morphism of v-sheaves 
        \[
            \coprod_{i\in I} \Spd(B_i,B_i^+) \twoheadrightarrow{}\mathcal{F}.
        \]
        \item[\textup{(2)}] \textbf{Formally separated:}
        The diagonal morphism $\Delta\colon\mathcal{F}\to \mathcal{F}\times \mathcal{F}$ is formally closed; that is, after any base change by a morphism from a formal Huber pair, the induced diagonal is a closed immersion.
    \end{enumerate}
    In this situation, there is a natural continuous \emph{specialization map}
    \[
        \mathrm{sp}_{\mathcal{F}}\colon |\mathcal{F}|\, \longrightarrow |\mathcal{F}_{\mathrm{red}}|,
    \]
    where $\mathcal{F}_{\mathrm{red}}$ denotes the reduced locus of $\mathcal{F}$.
\xqed{\blacktriangle}\end{dfn}

\begin{dfn}
    Let $\mathcal{F}$ be a specializing v-sheaf. Then $\mathcal{F}$ is called a \emph{pre-kimberlite} if:
    \begin{enumerate}
        \item[\textup{(1)}] The reduced locus $\mathcal{F}_{\mathrm{red}}$ is representable by a scheme.
        \item[\textup{(2)}] The natural adjunction morphism
        \[
            (\mathcal{F}_{\mathrm{red}})^{\diamond} \longrightarrow \mathcal{F}
        \]
        is a closed immersion.
    \end{enumerate}
\xqed{\blacktriangle}\end{dfn}

Let $\mathcal{F}$ be a pre-kimberlite v-sheaf, and let $C\subset |\mathcal{F}_{\textup{red}}|$ be a constructible subset. In \cite{gleason2020specialization}, Gleason defined the formal completion of $\mathcal{F}$ at $C$, denoted by $\mathcal{F}_{/C}\to \mathcal{F}$. Integral local Shimura varieties are examples of pre-kimberlite v-sheaves.

\begin{rmk}
    By the functoriality of specialization and completion, if $f\colon\CF\to \CF'$ is a map between kimberlite v-sheaves and $C\subset |\mathcal{F'}_{\textup{red}}|$ is a constructible subset, then $f$ induces an isomorphism from $\CF_{/f^{-1}(C)}$ to $\CF'_{/C}\times_{\CF'} \CF$, that we still denote by $f$.
\xqed{\lozenge}\end{rmk}

\begin{rmk}
    For a prismatic display \(\CD\), let \(\pi_{\textup{HT}}(\CD) \in \Gr_{(\CG,\mu)}\) denote the Hodge polygon of \(\CD\). By the deformation theory of displays \cite[Theorem 4.4.2]{ito2025deformationtheoryprismaticgdisplays}, there is an isomorphism \[\pi_{\textup{HT}} : \mathcal{M}^{\textup{int}}_{(\CG,\mu,b)/t} \xrightarrow{\;\sim\;} \Gr_{(\CG,\mu)/\pi_{\textup{HT}}(t)}.\]
\xqed{\lozenge}\end{rmk}

\begin{lem}\label{representabilityofcompletion}
    Let $\left(\CG,\mu,b\right)$ be a regularly sparse Shimura datum and $t\in |\mathcal{M}^{\textup{int}}_{(\CG,\mu,b),\,\textup{red}}|$. The map $\Flex\colon\mathcal{M}^{\textup{int}}_{(\CG,\mu,b)/t}\to \mathcal{M}^{\textup{int}}_{({\CG}^\triangleright_\rho,\mu^\triangleright,b)/\Flex(t)}$ is a closed immersion of formal schemes.
\end{lem}

\begin{proof}
    This map is a proper injective map between two affine formal schemes and, therefore, it is a closed immersion. 
\end{proof}

Now we want to prove the representability of $\Flex$. The proof is similar to the proof of \cite[Theorem 3.7.1]{pappas2022p}. Let $E$ be a local field, and let $\breve{E}$ denote the completion of its maximal unramified extension. 
We denote by $\breve{\mathcal{O}}$ the ring of integers of $\breve{E}$ and by $k$ its residue field. We begin by recalling the following lemma:

\begin{lem}{\cite[Proposition 3.7.3]{pappas2022p}}\label{criteriaforrepresentability}
    Assume that $\mathfrak{X}$ is a formal scheme over $\Spf(\breve{\mathcal{O}})$ and we have the data $\left(Z,T,\left(V_t\right)_{t\in T\left(k\right)}\right)$ where:
    \begin{enumerate}[label=(\roman*)]
        \item $Z\subset \mathfrak{X}^{\textup{rig}}$ is a closed immersion of rigid spaces.
        \item $T\subset \mathfrak{X}_{\textup{red}}$ is a reduced closed subscheme.
        \item $V_t\subset \mathfrak{X}_{/t}$ is a closed formal sub-scheme such that $sp\left(|Z|\right)\subset T\left(k\right)$ and $$\left(\mathfrak{X}_{/t}\right)^{\textup{rig}}\cap Z=V_t^{\textup{rig}}$$
    \end{enumerate} 
    Then there exists a unique closed formal subscheme $\mathfrak{Y}$ of $\mathfrak{X}$ such that $$\mathfrak{Y}^{\textup{rig}}=Z,\, \mathfrak{Y}_{\textup{red}}=T,\ \text{and}\,\, \mathfrak{Y}_{/t}=V_t$$
\end{lem}

\begin{thm}\cite[Theorem 3.3.3]{pappas2022p}\label{reducedlocus}
    The reduced locus $\left(\mathcal{M}^{\textup{int}}_{(\CG,\mu,b)}\right)_{\red}$ is representable by a perfect scheme $X_{(\CG,\mu,b)}$.  
\end{thm} 

\begin{rmk}
    The integral local Shimura datum $(\CG^\triangleright,\mu^\triangleright)$ is of Hodge type and therefore by \cite[Theorem 2.5.4]{pappas2022integral} the v-sheaf $\mathcal{M}^{\textup{int}}_{({\CG}^\triangleright_\rho,\mu^\triangleright,b)}$ is representable by the diamond of a smooth formal scheme.
\xqed{\lozenge}\end{rmk}

\begin{thm}\label{regularsparseisrepresentable2}
    Let $\left(\CG,\mu,b\right)$ be a regularly sparse Shimura datum. The v-sheaf $\mathcal{M}^{\textup{int}}_{(\CG,\mu,b)}$ is representable by the diamond of a smooth formal scheme.
\end{thm}

\begin{proof}
    By Lemmas \ref{colsedimmersionofvsheaf} and \ref{representabilityofcompletion}, we have closed immersions 
    \begin{enumerate}[label=(\roman*)]
        \item $\mathcal{M}_{(\CG,\mu,b)}\hookrightarrow  \mathcal{M}_{\{(\CG_\rho^\triangleright,\mu^\triangleright,b)\}}$
        \item $X_{(\CG,\mu,b)}\hookrightarrow X_{({\CG}^\triangleright_\rho,\mu^\triangleright,b)}$  
        \item $\mathcal{M}^{\textup{int}}_{(\CG,\mu,b)/t}\hookrightarrow  \mathcal{M}^{\textup{int}}_{({\CG}^\triangleright_\rho,\mu^\triangleright,b)/\Flex(t)}$
    \end{enumerate}

    Therefore, by Lemma \ref{criteriaforrepresentability}, there exists a formal scheme that represents $\mathcal{M}^{\textup{int}}_{(\CG,\mu,b)}$.
\end{proof}

\begin{thm}\label{regularsparseperfect2}
    Let $\left(\CG,\mu,b\right)$ be a regularly sparse Shimura datum. The local Shimura variety at infinite level $\mathcal{M}_{(\CG,\mu,b),\infty}\coloneqq \varprojlim\limits_{K}\, \mathcal{M}_{(\CG,\mu,b),K}$ is representable by a perfectoid space.
\end{thm}

\begin{proof}
    By Lemma \ref{colsedimmersionofvsheaf} we obtain a closed immersion $\Flex\colon\mathcal{M}_{(\CG,\mu,b)}\hookrightarrow \mathcal{M}_{(\CG_\rho^\triangleright,\mu^\triangleright,b)}$. By Corollary \ref{Flexarbirarylevel} we can lift this to a closed immersion $\Flex\colon\mathcal{M}_{(\CG,\mu,b),\infty}\hookrightarrow \mathcal{M}_{(\CG_\rho^\triangleright,\mu^\triangleright,b),\infty}$. Since the local Shimura variety $\mathcal{M}_{(\CG_\rho^\triangleright,\mu^\triangleright,b)}$ is of Hodge type, it follows that it is representable by a perfectoid space \cite{scholze2015torsion}. Therefore, $\mathcal{M}_{(\CG,\mu,b),\infty}$ is also representable by a perfectoid space.
\end{proof}

\section{Integral Model of Meta-Unitary Shimura Varieties}\label{integralmodelsection}

Let $(G,\mu)$ be a meta-unitary Shimura datum in the sense of Definition \ref{metaunitarydatumdfn}. 
In \cite{bultel2008pel}, B\"ultel provides a modular description of the special fiber of the would-be integral model of $\Sh(G,\mu)$, and then uses deformation theory to lift the special fiber and construct the integral model. 
However, he does not provide an explicit description of the integral model itself; for a summary and the notations used in this section, please see Appendix \ref{defmeta-unitarysection}. 
Let $S$ denote the integral model constructed in \cite{bultel2008pel}, and $\CS$ be the associated formal scheme. 
Our goal is to give an explicit description of $\mathcal{S}^\diamond$, analogous to the description of the special fiber. 
To this end, we first define a perfectoid version of the PEL-type moduli problem considered by B\"ultel.

Let $\CG$ be a reductive model for $G$ over $\BZ_p$, and $K^p\subset G(\mathbb{A}^{\infty,p})$ a compact open subgroup. Let $K=K^p\CG(\BZ_p)$ be a full level.

\begin{dfn}\label{DefModuliProb}
    Let $L$ be a $\textup{CM}$-field and $\mathcal{O}_L$ its ring of integers. Let $R_1,R_2,R_3$ and $R_{123}$ be involutive algebras that are free and of finite rank over $\CO_L$. Consider the moduli problem $\CS_{K^p}^{\textup{tan},\diamond}$ over connected affinoid perfectoid pairs that sends $\Spa(R,R^+)$ to the following set:
    
    \begin{itemize}
        \item A quadruple $\left(Y,\lambda,\eta,\iota\right)$ consisting of a polarized Abelian variety $(Y,\lambda)$ on $R^+$ with $K^p$-level structure $\eta$ and an $\mathcal{O}_L^{\oplus3}$-action $\iota$ that satisfies the usual conditions.
        
        \item Actions $y_i$ of $R_i$ on $\bar{Y}_i=(e_i.Y)\times_{\Spec R^+} \Spec (R^+/p)$ that satisfy the determinant condition (where $e_i\in \mathcal{O}_L^{\oplus3}$ are the standard basis elements), and an action of $R_{123}$ on $\bar{Y}_{123}=e_1\bar{Y}_1\otimes e_2\bar{Y}_2\otimes e_3\bar{Y}_3$ (please check the construction \ref{hypothesisconstruction} for the definition of the tensor product) such that the action of  $R_{123}$ commutes with its action on $V_{123}\otimes \mathbb{A}^{\infty,p}$. More precisely, for $s_0\in S$, the composition:
    
        \[
            \begin{tikzcd}
            \mathbb{A}^{\infty,p}\otimes V_{123} \arrow{r}{\eta}& \bigotimes H_{\text{ét}}^1\left(Y_{i,s_0},\mathbb{A}^{\infty,p}\right)\arrow{r}& H_{\text{ét}}^1\left(Y_{123,s_0},\mathbb{A}^{\infty,p}\right) 
            \end{tikzcd}
        \]
        commutes with the action of $R_{123}$.
    \end{itemize}
\xqed{\blacktriangle}\end{dfn}

\begin{rmk}
    This moduli problem is a closed subsheaf of the integral model of the Shimura variety of polarized Abelian varieties with an action of $\mathcal{O}_L^{\oplus3}$. Therefore we obtain an embedding $\CS_{K}^{\textup{tan},\diamond}\to \Sh_K^\diamond(G^\triangleright,\mu^\triangleright)$
\xqed{\lozenge}\end{rmk}

We also need to define the functor $\Flex$ over a formal smooth scheme. We will denote by $\CB^{\textup{Witt, nil}}_{(\CG,\mu)}$ the substack of $\CB^{\textup{Witt}}_{(\CG,\mu)}$ consisting of adjoint-nilpotent displays  (see \cite[Definition 4.3]{bueltel2020g} for the definition of an adjoint-nilpotent $(\CG,\mu)$-display).

\begin{lem}
    Let $(\CG,\mu)$ be a sparse integral Shimura datum. Then the functor $\Flex$ extends to a functor \[\Flex \colon \CB^{\textup{Witt, nil}}_{(\CG,\mu)} \longrightarrow \CB^{\textup{Witt, nil}}_{(\CG_\rho^\triangleright,\mu^\triangleright)} \] when both stacks are restricted to smooth $p$-adic formal schemes. 
\end{lem}

\begin{proof}
    Let $\FX$ be a smooth $p$-adic formal scheme. We may cover $\FX$ by affine smooth formal schemes $U_i = \Spf B_i$ such that each $B_i$ admits a pro-\'etale cover by a perfectoid ring $B_i^\infty$ \cite[§2.3]{takaya2024moduli}.  

    Let $\CD$ be an adjoint-nilpotent prismatic $(\CG,\mu)$ Witt display over $B_i$. After base change to $B_i^\infty$, we obtain an adjoint-nilpotent prismatic display $\Flex(\CD_{B_i^\infty})$ over $B_i^\infty$, and hence, by \cite{bartling2022mathcal}, a Witt display over $B_i^\infty$.
    
    The display $\Flex(\CD_{B_i^\infty})$ satisfies descent with respect to the pro-\'etale cover $B_i^\infty/B_i$. Therefore, it descends to a Witt display $\Flex(\CD_{i})$ over $B_i$. Finally, since these local constructions are compatible on overlaps, the displays $\Flex(\CD_{i})$ glue to yield a global Witt display $\Flex(\CD)$ over $\FX$.
\end{proof}

\begin{dfn}\label{hyperspecial}
    Let $\mathcal{S}^{\textup{mixed}}$ be the v-sheaf making the following diagram Cartesian:
        \begin{equation}\label{DiagMixedShtuka}
\begin{tikzcd}
\mathcal{S}_{K}^{\mathrm{mixed}} \arrow[r, "\pi"] \arrow[d] & \CB^{\textup{perf}}_{(\CG,\mu)} \arrow[d, "\Flex"] \\
\mathcal{S}^{\mathrm{tan},\diamond}_{K^p} \arrow[r] & \CB_{(\CG^\triangleright_\rho,\mu^\triangleright)}^{\mathrm{tan}}
\end{tikzcd}
\end{equation}

\xqed{\blacktriangle}\end{dfn}

\begin{rmk}
    By definition, we have an isomorphism $\mathcal{S}^{\textup{mixed}}_{\bar{\BF}_p}\cong\mathcal{S}^\diamond_{\bar{\BF}_p}$.
\xqed{\lozenge}\end{rmk}

\begin{thm}\label{diamond}
    The v-sheaves $\mathcal{S}_{K}^\diamond$ and $\mathcal{S}_{K}^{\textup{mixed}}$ are naturally isomorphic. 
\end{thm}

\begin{proof}
  Let $\mathcal{D}$ be the universal $(G,\mu)$-display over $\mathcal{S}_{K}$, and let $\bar{\mathcal{A}}$ denote the universal Abelian variety over $\mathcal{S}_{K^p,\BF_q}$. The display $\Flex(\mathcal{D})$ lifts the display associated with $\bar{\mathcal{A}}$. By the Serre–Tate and Grothendieck–Messing theorems, there exists a unique lift of $\bar{\mathcal{A}}$ whose display is $\Flex(\mathcal{D})$. Thus, we obtain a morphism $\mathcal{S} \to \mathcal{S}^{\mathrm{tan},\diamond}_{K}$ that makes the following diagram Cartesian:
    \[
    \begin{tikzcd}
    \CS^\diamond_{K} \arrow{r} \arrow{d} & \CB^{\textup{perf}}_{(\CG,\mu)} \arrow[d, "\Flex^{\mathrm{tan}}"] \\
    \mathcal{S}^{\mathrm{tan},\diamond}_{K} \arrow{r} & \CB_{(\CG_\rho^\triangleright,\mu^\triangleright)}^{\mathrm{tan}}
    \end{tikzcd}
    \]
    Therefore, $\CS_{K}^{\diamond}$ is isomorphic to $\CS_{K}^{\textup{mixed}}$.
\end{proof}

Recall the definition of $G_\rho^\triangleright$ from \ref{dfnGtriangle}.
\begin{cor}\label{closed immersion2}
    There exist compact open subgroups $K_p \subset G(\BQ_p)$, $\tilde{K} \subset G_\rho^\triangleright(\BA^{\infty})$, and a Zariski closed immersion
    \[
    \Flex \colon \Sh_{K^pK_p}(G,\mu) \xhookrightarrow{\quad} \Sh_{\tilde{K}}(G^\triangleright,\mu^\triangleright).
    \]
\end{cor}

\begin{proof}
    By Theorems \ref{diamond} and \ref{complex_shimura_variety}, there is a morphism
    \[
    \Sh_{K^pK_p}(G,\mu)^\diamond \longrightarrow \Sh_{\tilde{K}}(G_\rho^\triangleright,\mu^\triangleright)^\diamond
    \]
    which is obtained by base change from $\Flex^{\textup{tan}}$. Since $\Flex^{\textup{tan}}$ is a monomorphism, and a monomorphism between proper diamonds is a closed immersion, this map is a closed immersion.

    For a scheme $X$, we denote its rigidification by $X^{\textup{rig}}$. Both $\Sh_{K^pK_p}(G,\mu)^{\textup{rig}}$ and $\Sh_{\tilde{K}}(G_\rho^\triangleright,\mu^\triangleright)^{\textup{rig}}$ are normal rigid spaces. By \cite{scholze2020berkeley}, the functor from normal rigid spaces to diamonds is fully faithful. Hence, we obtain a closed immersion
    \[
    \Flex \colon \Sh_{K^pK_p}(G,\mu)^{\textup{rig}} \xhookrightarrow{\quad} \Sh_{\tilde{K}}(G_\rho^\triangleright,\mu^\triangleright)^{\textup{rig}}.
    \]

    Since the rigid space $\Sh_{\tilde{K}}(G_\rho^\triangleright,\mu^\triangleright)^{\textup{rig}}$ is projective, GAGA implies that $\Flex$ arises from a closed immersion between the corresponding Shimura varieties.
\end{proof}

\section{Meta-unitary Shimura Varieties at Deeper Levels}\label{deeperlevelsection}

If we want to understand the geometry of a Shimura variety at the $p^\infty$-level, we have to think about levels that are neither hyperspecial nor even parahoric. At these levels, we do not have a good theory of integral models for Shimura varieties. However, for every level $K$, we still have a shtuka with $K$-level structure over the Shimura variety with $K$-level structure. Let $\CG$ be a reductive group over $\BZ_p$.

\begin{dfn}
    Let $\CF$ be a v-sheaf over $\Spd(\BZ_p)$. A \emph{$\CG$-shtuka over $\CF$ with one leg} is a rule assigning to every perfectoid space $S \to \CF$ a $\CG$-shtuka with one leg over $S$. These data must be compatible with pullback in $S$. Equivalently, a $\CG$-shtuka over $\CF$ is a map from $\CF$ to the v-stack of $\CG$-shtukas.
\xqed{\blacktriangle}\end{dfn}

\begin{dfn}
    Let $\FX$ be a scheme over $\mathbb{Q}_p$ or a $p$-adic formal scheme over $\Spf(\BZ_p)$, and let $\FX^\diamond$ be its associated diamond over $\Spd(\BZ_p)$. A \emph{$\CG$-shtuka over $\FX$} is a $\CG$-shtuka over $\FX^\diamond$.
\xqed{\blacktriangle}\end{dfn}

\begin{rmk}
    A $\CG$-shtuka over a scheme or formal scheme induces a $\CG(\BZ_p)$-local system on the associated diamond $S^\diamond$. If $S$ is an analytic adic space (for example if it is a scheme over $\BQ_p$), then, by the equivalence of the \'etale sites of $S$ and $S^\diamond$, this $\CG(\BZ_p)$-local system corresponds to a $\CG(\BZ_p)$-local system on $S$ itself.
\xqed{\lozenge}\end{rmk}

\begin{dfn}\label{shtukaoverformalschemedfn}
    Let $K\subset \CG(\BZ_p)$ be an open compact subgroup. A $\CG$-shtuka with $K$-level structure over a $\BQ_p$-scheme $S$ consists of a triple $(\SP,T,\iota)$, where $\SP$ and $T$ are respectively a $\CG$-shtuka and a $K$-torsor over $S$, and $\iota$ is an isomorphism between the induced $\CG(\BZ_p)$-torsor $T \times^{K} \CG(\BZ_p)$ and the $\CG(\BZ_p)$-torsor associated with $\SP$.
\xqed{\blacktriangle}\end{dfn}

\begin{dfn}
  A de Rham $\CG$-torsor is a $\CG$-torsor $\mathcal{P}$ such that for every representation of $\CG$, the associated local system is de Rham. 
\xqed{\blacktriangle}\end{dfn}

\begin{rmk}
    Using Beauville--Laszlo theorem, one can associate a $\CG$-shtuka to every de Rham $\CG$-torsor. For details, see \cite[Definitions 2.6.5, 2.6.6]{pappas2022p}.
\xqed{\lozenge}\end{rmk}

\begin{thm}{\cite{Hansen2016PeriodMA}}
    Fix open compact subgroups $K^p \subset G(\mathbb{A}^{\infty,p})$ and $K_p \subset G(\mathbb{Q}_p)$. Consider the associated Shimura variety $\Sh_{K_pK^p}(G,\mu)$ and its pro-\'etale cover \[\Sh_{K_{p,\infty}K^p}(G,\mu) \coloneqq \varprojlim_{K'_p \subset K_p} \Sh_{K'_pK^p}(G,\mu).\] Then the $K_p$-torsor naturally associated with $\Sh_{K_{p,\infty}K^p}$ is de Rham, and there exists a canonical period map \[\pi_{\mathrm{dR}} : \Sh_{K_{p,\infty}K^p}(G,\mu) \longrightarrow \Gr_{(G,\mu)}.\]
\end{thm}

\begin{rmk}\label{RemSHtoSHT}
    Let $K^p$ be  a fixed level away from $p$, as observed in \cite{pappas2022p}, for every level $K=K_p$, the de Rham period morphism
    \[
    \pi_{\mathrm{dR}} \colon \Sh_{K^pK_{p^\infty}}(G,\mu) \longrightarrow \Gr_{(G,\mu)}
    \]
    yields a $\CG$-shtuka with $K_p$-level structure on $\Sh_{K_pK^p}(G,\mu)$. We denote the resulting morphism of diamonds by
    \[
    \pi_{\textup{crys}} \colon \Sh_{K_pK^p}(G,\mu)^\diamond \longrightarrow \Sht_{(\CG,\mu),K_p}.
    \]
\xqed{\lozenge}\end{rmk} 

\begin{assumption*}[CS]
    The morphism $\pi_{\textup{crys}}$ and the morphism $\Fib_G\circ\pi$ in the adic generic fiber of the diagram (\ref{DiagMixedShtuka}) are the same (here $\Fib_G$ is the fully faithful functor from the stack of $G$-displays to $G$-shtukas, cf. Construction \ref{ConsFib}). 
\end{assumption*}

\begin{rmk}
    In a work in progress \cite{hedayatzadehrappaportlanglands}, we prove that the Assumption (CS) always holds.
\xqed{\lozenge}\end{rmk}

\begin{pro}\label{deeperlevel}
   Let $K = K_pK^p$ and $K' = K'_pK^p$ be open compact subgroups of $G(\BA^\infty)$ with $K_p \subset K'_p \subset G(\BQ_p)$. Then the following diagram of diamonds is Cartesian:
        \[
        \begin{tikzcd}
        \Sh_{K'}(G,\mu)^\diamond \arrow[r, "\pi_{\textup{crys}}"] \arrow[d] & \Sht_{(\CG,\mu),K'} \arrow[d] \\
        \Sh_{K}(G,\mu)^\diamond \arrow[r, "\pi_{\textup{crys}}"'] & \Sht_{(\CG,\mu),K}
        \end{tikzcd}
        \] 
\end{pro}

\begin{proof}
    Note that the diagram is commutative by the compatibility of the associated local systems under change of level. Set 
    \[
        X \coloneqq \Sh_{K}(G,\mu)^\diamond \times_{\Sht_{(\CG,\mu),K}} \Sht_{(\CG,\mu),K'}.
    \] 
    Commutativity of the diagram yields a map $f\colon \Sh_{K'}(G,\mu)^\diamond \to X$. We want to prove that it is an isomorphism. We first show that it is surjective. Note that the map $\Sh_{K'}(G,\mu)\to \Sh_{K}(G,\mu)$ is surjective, and its fibers are in natural bijection with the cosets $K'/K$. By definition, the $G$-Shtukas over these Shimura varieties arise from the local systems associated with the covers
    \[
        \Sh_{K_{p,\infty}K^p}(G,\mu)\to\Sh_K(G,\mu) \quad \text{and} \quad \Sh_{K'_{p,\infty}K^p}(G,\mu)\to\Sh_{K'}(G,\mu)
    \]
    note that the two towers $\Sh_{K_{p,\infty}K^p}(G,\mu)$ and $\Sh_{K'_{p,\infty}K^p}(G,\mu)$ are naturally isomorphic. For a geometric point $x$ of $\Sh_{K'}(G,\mu)$, the fibers of the natural map between these local systems are again in bijection with $K/K'$; therefore, the map $f$ is surjective.
    
    Both $\Sh_{K'}(G,\mu)^\diamond$ and $X$ are finite \'etale covers of 
    $\Sh_{K}(G,\mu)^\diamond$ with the same number of sheets, namely $|K'/K|$.  
    It follows that this map is an isomorphism; therefore, the diagram is Cartesian.
\end{proof}

Our goal is to have a description of meta-unitary Shimura varieties at a cofinal subset of levels, similar to the description of hyperspecial levels:

\begin{dfn}
    Assume that we have a regular meta-unitary Shimura datum for $\left(G,\mu\right)$, and $\{K_{p,i}\},i\in {\{1,2,3,123\}}$ are compact open subgroups of $\textup{GU}\left(V_i,\Psi_i\right)(\mathbb{Q}_p)$. A good level at $p$ is a compact open subgroup $K_p$ obtained by taking the intersection of the inverse image of these level subgroups in $G$.
\xqed{\blacktriangle}\end{dfn}

In what follows, we will only work with good level structures. By taking the generic fiber of the diagram (\ref{DiagMixedShtuka}), we obtain the Cartesian diagram:
 \[
    \begin{tikzcd}
\Sh_{K}\left(G,\mu\right)^\diamond \arrow{r} \arrow{d} & \CB^{\textup{perf}}_{(\CG,\mu),\eta} \arrow{d} \\
\mathcal{S}^{\mathrm{tan},\diamond}_{K} \arrow{r} & \CB_{(\CG_\rho^\triangleright,\mu^\triangleright),\eta}^{\textup{tan}}
    \end{tikzcd}
    \]

\begin{thm}
    Let $\left(G,\mu\right)$ be a regular meta-unitary Shimura variety, and $K'_p\subset G(\BQ_p)$ a compact open subgroup. Set $K'=K'_pK^p$. Under Assumption (CS), there is a Cartesian diagram:
    
    \[
    \begin{tikzcd}
    \Sh\left(G,\mu\right)_{K'}^\diamond\arrow{r}\arrow{d}& \Sht_{(\CG,\mu),K'_p}\arrow{d}{\Flex}\\
    \mathcal{S}^{\mathrm{tan},\diamond}_{K'}\arrow{r}  & \Sht_{(\CG_\rho^\triangleright,\mu^\triangleright),K'_p}^{\textup{tan}}
    \end{tikzcd}
    \]
\end{thm}

\begin{proof}
    By Proposition \ref{deeperlevel} and Definition \ref{hyperspecial}, we have Cartesian diagrams
    \[
    \begin{tikzcd}
        \Sh_{K'}\left(G,\mu\right)^\diamond \arrow{r}\arrow{d}&\Sht_{(\CG,\mu),K'_p}\arrow{d}\\
        \Sh_K\left(G,\mu\right)^\diamond \arrow{r}&\Sht_{(\CG,\mu),K_p}\\
    \end{tikzcd}
    \]
    \[
    \begin{tikzcd}
        \mathcal{S}^{\mathrm{tan},\diamond}_{K'}\arrow{r}\arrow{d}  & \Sht_{(\CG^\triangleright,\mu^\triangleright),K'_p}^{\textup{tan}}\arrow{d}\\
        \mathcal{S}^{\mathrm{tan},\diamond}_{K}\arrow{r}  & \Sht_{(\CG^\triangleright,\mu^\triangleright)}^{\textup{tan}}\\
    \end{tikzcd}
    \]
    \[
\begin{tikzcd}
\Sh_{K}\left(G,\mu\right)^\diamond \arrow{r} \arrow{d} & \Sht_{(\CG, \mu)} \arrow{d} \\
\mathcal{S}^{\mathrm{tan},\diamond}_{K} \arrow{r} & \Sht_{(\CG_\rho^\triangleright,\mu^\triangleright)}^{\textup{tan}}
\end{tikzcd}
\]

Combine these with the cartesian diagram:
\[
\begin{tikzcd}
\CB^{\textup{perf}}_{(\CG,\mu),K'} \arrow{r} \arrow{d} & \Sht_{(\CG, \mu,K')} \arrow{d} \\
\mathcal{B}^{\textup{perf}}_{(\CG,\mu),K} \arrow{r} & \Sht_{(\CG,\mu),K}
\end{tikzcd}
\]

\end{proof}

\begin{thm}\label{globalperfectoid2}
    Let $\left(G,\mu\right)$ be a regular meta-unitary Shimura datum and $K^p$ be a level away from $p$. Under the assumption (CS), the Shimura variety at infinite level $$\Sh_{K^p}\left(G,\mu\right)\coloneqq \varprojlim_{K_p} \Sh_{K_pK^p}\left(G,\mu\right)$$ is representable by a perfectoid space.
\end{thm}

\begin{proof}
    The map $\Flex$ is a closed immersion to $\CS_{K}$, and there is an embedding of this variety to the Siegel modular variety by this variety. We have a closed immersion from $\Sh\left(G,\mu\right)_{K^p}$ to the Siegel modular variety at infinite level that is perfectoid by \cite[Theorem 1.6]{scholze2015torsion}.
\end{proof}


\begin{appendices}

\section{Prismatic Display Groups}\label{PrismDispGrp}

In this section, we briefly recall the construction of the prismatic display group defined in \cite{hedayatzadeh2024deformations} and indicate how to extend it to the setting of multi-leg prisms.

Let $(A, I)$ be an $\BZ_q$-prism. The \emph{Nygaard filtration} on $(A, I)$ is defined by
\[
\CN^i_A \coloneqq \sigma^{-1}(I^i A),
\]
where $\sigma$ denotes the Frobenius lift on $A$.

Let $\CG$ be a smooth affine group scheme over $\mathbb{Z}_p$, and let $\mu \colon \BG_{m} \to \CG_{\BZ_q}$ be a cocharacter defined over $\BZ_q$. This action defines a decreasing filtration $\CF^\bullet $ on $\CO_{\CG}\otimes_{\BZ_p}\BZ_q$ by
\[
\CF^j \coloneqq \bigoplus_{i \ge j} \CO_{G,i},
\]
where $\CO_{\CG}\otimes_{\BZ_p}\BZ_q= \bigoplus_{i \in \mathbb{Z}} \CO_{\CG,i}$ denotes the weight decomposition of $\CO_{\CG}\otimes_{\BZ_p}\BZ_q$ via $\mu$.

\begin{dfn}
    We define the \emph{prismatic display group} $\CG^\mu(A,I)$ to be subgroup of $\CG(A)=\Hom(\CO_{\CG}\otimes_{\BZ_p}\BZ_q,A)$ consisting of those morphisms $k \colon \CO_{\CG}\otimes_{\BZ_p}\BZ_q \to A$ which respect the filtrations on both sides, that is,
    \[
    k(\CF^j) \subseteq \CN^j_A, \quad \text{for all } j \ge 0.
    \]
    For a multi-leg prism $(A,\mathbf{I})$, and a cocharacter $\boldsymbol{\mu}$ of $\prod_{i\in\BZ/r}\CG$, we define the associated \emph{prismatic display group} to be \[\CG^{\boldsymbol{\mu}}(A,\mathbf{I})\coloneqq \prod_{i\in\BZ/r}\CG^\mu(A,I_i)\]
\xqed{\blacktriangle}\end{dfn}

\begin{dfn}
    Assume that $I$ is generated by $d\in A$. The \emph{display map}
    \[
    \Phi_d^\mu \colon \CG^\mu(A,I) \longrightarrow \CG(A[1/d])
    \]
    is defined by the formula
    \[
    \Phi_d^\mu(k) \coloneqq \mu(d) \, \sigma(k) \, \mu(d^{-1}).
    \]
    It is straightforward to verify (see \cite[Appendix A]{hedayatzadeh2024deformations}) that $\Phi_d^\mu(k)$ in fact lies in $\CG(A)$, giving rise to a well-defined morphism
    \[
    \Phi_d^\mu \colon \CG^\mu(A,I) \to \CG(A).
    \]
\xqed{\blacktriangle}\end{dfn}

Now, we want to show that we can make the above constructions independent from the choice of orientation. Following \cite[Definition 4.5]{hedayatzadeh2024deformations}, this is achieved by taking an inverse limit over all possible generators.

\begin{construction}
    For a principal ideal $I \subset A$, let $d$ range over all generators of $I$. For each such $d$, set $\CG(A)_d \coloneqq \CG(A)$. Whenever $e = vd$ with $v \in A^\times$, define a transition morphism
    \begin{align*}
       \alpha_{e/d}\colon \CG(A)_d \to &\, \CG(A)_{e}\\
       X \mapsto &\, X\mu(v^{-1})
    \end{align*}   
    The \emph{canonical limit group} associated to $(A,I)$ is then defined as
    \[\CG(A, I) \coloneqq \varprojlim_d\, \CG(A)_d\]
    where the limit is on the set of generators of $I$ and the transition morphisms are the $\alpha$ defined above.

    For each generator $d$ of $I$, we have a \emph{$\Phi_d^\mu$-conjugation} on $\CG(A)_d$: for $k\in \CG^{\mu}(A,I)$ and $X\in \CG(A)_d$, we set $k\cdot X\coloneqq k^{-1}X\Phi^{\mu}_d(k)$. One can see that the transition morphisms $\alpha_{e/d}$ are $\CG^\mu(A,I)$-equivariant, and so, we obtain an action of $\CG^\mu(A,I)$ on $\CG(A,I)$, that we will call \emph{$\Phi^\mu_I$-conjugation}.

    Now, for a multi-leg prism $(A,\mathbf{I})$, we define the \emph{canonical limit group} \[\CG(A, \mathbf{I}) \coloneqq \prod_{i \in \mathbb{Z}/r\mathbb{Z}} \CG(A, I_i).\]
    The $\Phi^\mu_I$-conjugations yield an action of $\CG^{\boldsymbol{\mu}}(A,\mathbf{I})$ on $\CG(A,\mathbf{I})$, component-wise, and we call it $\Phi^{\boldsymbol{\mu}}_{\mathbf{I}}$-conjugation.
\xqed{\blacktriangledown}\end{construction}

\begin{rmk}
    The group functors $\CG$ and $\CG^\mu$ are sheaves on the perfectoid site.
\xqed{\lozenge}\end{rmk}

\begin{rem}
    For $r = 1$, the above definitions specialize to those of \cite{hedayatzadeh2024deformations}. The inverse limit over generators replaces the ad hoc dependence on a single element $d$, making the construction functorial in the underlying prism.
\end{rem}

\section{Meta-unitary Shimura Varieties}\label{defmeta-unitarysection}

The main motivation of \cite{bultel2008pel} for the definition of the functor $\Flex$ in equal characteristic is to construct the special fiber of a certain class of Shimura varieties called \emph{meta-unitary} Shimura varieties. This class contains some Shimura varieties that are not of Abelian type, including the ones in Example \ref{example}. 

In this section, we recall the definition of meta-unitary Shimura varieties and the construction of an integral model. There are many conditions involved in this definition; most of them are technical, and one can hope that they can be removed. The important thing is that for the display datum of this class of Shimura varieties, we can define the functor $\Flex$.

In the constructions and definitions that follow, we will refer to the following two conditions. Let $W_1$, $W_2$, and $W_3$ be vector spaces, and $\boldsymbol{\mu}_i=(\mu_{i,j})_{j=1}^r$ a cocharacter of $\GL(W_i)^r$ ($i=1,\dots,3$). Consider the following two conditions:

\begin{cond}\label{hypothesistensor}
    For each $j$ there exists $i\in\{1,2,3\}$ such that the $j$-th component of $\boldsymbol{\mu}_i$ is  either the trivial cocharacter or is equal to $z\mapsto \left(z,z,\ldots,z\right)$.
\end{cond}

\begin{cond}\label{hypothesiscompactness}
    There exists a component $j$ such that, for every $i \in \{1,2,3\}$, the $j$-th component of $\boldsymbol{\mu}_i$ is either the trivial cocharacter or the cocharacter $z \mapsto (z, z, \ldots, z)$.
\end{cond}

First, let us recall an important class of sparse local Shimura data defined in \cite{bultel2008pel}. Let $\CG_0$ be a smooth group scheme defined over $\BZ_{p^r}$, and $\boldsymbol\mu$ be a cocharacter of $\CG_0^r$. Recall the notations from \ref{sectionsparse}. We are interested in the case where $\CG_0$ is the stabilizer of the action of some involutive algebras.

\begin{dfn}\label{gaugeddisplaydatum}
    A \emph{gauged $\left(\CG_0,\mu\right)$-display datum} consists of $\CG_0$-representations $\left(V_1,\rho_1\right),\left(V_2,\rho_2\right)$ and $\left(V_3,\rho_3\right)$ in free $\mathbb{Z}_p$-modules of finite rank satisfying the following conditions, where we denote by $\left(V_{123},\rho_{123}\right)$ the tensor product of $(V_{i},\rho_{i})$ and $\mu_i=\rho_i\circ\mu$:
        \begin{enumerate}
            \item For each $i\in\{1,2,3\}$ the display datum $\left(\GL\left(V_i\right)^\Gamma,\mu_i^\Gamma\right)$ is sparse.
            
            \item There exist involutive algebras $R_1,R_2,R_3$ and $R_{123}$ that are free and of finite rank over $\mathbb{Z}_{p^r}$ and an action of $R_j$ on $V_j$ for $j\in\left\{1,2,3,123\right\}$ such that $G_0$ is the  stabilizer of these actions inside $\prod \GL(V_j)$.            
            \item If $\left(\GL\left(V_i\right)^r,\boldsymbol{\mu}^\triangleright_i\right)$, for $i\in\{1,2,3\}$, is the minuscule display datum given by Construction \ref{tildemu} applied on $\left(\CG_0,\mu\right)$, then $\left(\GL\left(V_{123}\right)^r,\boldsymbol{\mu}^\triangleright_{1}\otimes\boldsymbol{\mu}^\triangleright_2\otimes\boldsymbol{\mu}^\triangleright_3\right)$ is minuscule. Equivalently, the cocharacters $\boldsymbol{\mu}^\triangleright_{i}$ satisfy Condition \ref{hypothesistensor}.
            \item Characters $\boldsymbol{\mu}^\triangleright_{i}$ satisfy Condition \ref{hypothesiscompactness}. 
        \end{enumerate}
\xqed{\blacktriangle}\end{dfn}

Let us comment on why one needs such conditions:

\begin{rmk}
    Condition $1$ is the most important condition because it allows us to define the map $\Flex$. The purpose of considering three representations and their tensor product satisfying Conditions $2$ and $3$ instead of a single representation is that we want to remain in a PEL-type setting and still have interesting examples satisfying Condition $1$. The last condition is also a useful technical one ensuring that all the $(\Res_{\mathbb{Z}_{p^r}/\mathbb{Z_p}}\GL_n,\boldsymbol\mu^\triangleright)$-displays are adjoint-nilpotent.
\xqed{\lozenge}\end{rmk}

There is also a skew-Hermitian variant of the notion of gauged display datum and the functor $\Flex$, which is more useful for global applications.
    
\begin{dfn}\label{skewhermitiandisplaydatm}
    Let $R$ be a finite \'etale $\BZ_q$-algebra of rank $2$. A \emph{skew-hermitian gauged display datum over $R$} is the data of three unitary representations of $\CG_0$ over finite free $R$-algebras that satisfy the above conditions.
\xqed{\blacktriangle}\end{dfn}

\begin{rmk}\label{skewhermitianvarianofFLex}
    One can define a skew-Hermitian variant of the functor $\Flex$: Using the Tannakian point of view, a \emph{skew-Hermitian display} is a triple consisting of a display $\CD$, a display $\CK$ for the multiplicative group $\BG_{m}$, and a skew-Hermitian pairing on $\CD$. We then apply the functor $\Flex$ to $\CD$ and $\CK$ and transfer the skew-Hermitian structure to $\Flex(\CD)$ using the compatibility of $\Flex$ with duality and pairings. For details, see \cite[Section 6.5]{bultel2008pel}.
\xqed{\lozenge}\end{rmk}

\begin{rmk}\label{RemSht^R}
    In \cite{bultel2008pel}, B\"ultel denote $\Flex^{\textup{tan}}$ and $\Sht^{\textup{tan}}$ by $\Flex^R$ and $\Sht^R$ respectively.
\xqed{\lozenge}\end{rmk}

\begin{setup}
    Let $L$ be a CM-field and $L^+$ be its maximal totally real subfield. For simplicity, let us assume that $L^+/\mathbb{Q}$ is Galois, although it is not necessary. Denote the set of embeddings of $L$ into $\mathbb{C}$ by $L^{\textup{an}}$. Let $\alpha\mapsto \bar{\alpha}$ be the generator of $\textup{Gal}\left(L/L^+\right)$. Let $V$ be a vector space over $L^+$ together with a skew-Hermitian form $\Psi$ on $V\otimes_{L^+}L$. Let $\textup{GU}(V,\psi)$ be the group scheme over $L^+$ with the functor of points: 

 \begin{align*}
    \textup{GU}(V,\psi)(R) = \{\, g \in \GL(V_{R\otimes_{L^+}L}) \mid {}& \psi(gv, gw) = \lambda(g) \, \psi(v,w) \\
    &\text{ for some } \lambda(g) \in F^\times,\ \forall v,w \in V \,\}.
\end{align*}
    
    A Hodge structure of weight $-1$ is said to be \emph{skew-Hermitian with coefficients in $L$} if the map $\mathbb{S}\to \GL(\mathbb{V}_{\mathbb{R}})$ factors through $\textup{GU}\left(V,\Psi\right)_{\mathbb{R}}$.    
\end{setup}

A Hodge structure with coefficients in $L$ gives us a graded Hodge structure: $$V_{\mathbb{C}}=\bigoplus_{\tau\in L^{\textup{an}}}\bigoplus_{p+q=-1} V_\tau^{p,q}$$ where $\overline{V_\tau^{p,q}}={V}_{\bar{\tau}}^{q,p}$.

Let $\sigma$ be a Frobenius lift on $L$, and $S$ an orbit of $L^{an}$ under the action of $\sigma$. Choose a bijection $j_S\colon\mathbb{Z}/|S|\mathbb{Z}\to S$, and consider the double graded module $$V_{j_S}\coloneqq\bigoplus_{i\in\frac{\BZ}{|S|\BZ}}\bigoplus_{p+q=-1} V_{j_S\left(i\right)}^{p,q}$$ 

\begin{dfn}
    Let $\boldsymbol\mu_S$ be the cocharacter obtained by the above decomposition on $V_S$. Let $\Gamma\subset S$ be the subset consisting of $\gamma\in S$ such that $\mu_{\gamma}\not=\id$. From now on we use the notations of \ref{notationgamma}. We say that $j_S$ is sparse if $\mu_S$ is sparse in the sense of Definition \ref{notationgamma}.  
\xqed{\blacktriangle}\end{dfn}

 We say that $V$ is \emph{sparse} if, for each orbit $S$, there is a sparse bijection $j_S$ such that the bijections $j_S$ and $j_{\bar{S}}$ are compatible with complex conjugation:

    \begin{itemize}
        \item If $S \not= \bar{S}$, then $\overline{j_S\left(i\right)} = {j_{\bar{S}}\left(i\right)}$.
        \item If $S = \bar{S}$, then $\overline{j_S\left(i\right)} = {j_{S}\left(|S|-i\right)}$.
    \end{itemize}
    
\begin{construction}\label{archemedianmutilde}
    Let $(V,h,\Psi)$ be a sparse skew-Hermitian Hodge structure with coefficients in $L$. One can associate with $V$ a minuscule Hodge structure $(V^\triangleright,h^\triangleright)$ using Construction \ref{tildemu}. Furthermore, one can upgrade $V^\triangleright$ to a skew-Hermitian Hodge structure $(V^\triangleright,h^\triangleright,\Psi^\triangleright)$ such that there exists an isomorphism $\epsilon^{\infty,p}\colon(V\otimes \BA^{\infty},\Psi)\to (V^\triangleright\otimes \BA^{\infty},\Psi^\triangleright)$. For more details see \cite[Lemma 8.2]{bultel2008pel}.
\xqed{\blacktriangledown}\end{construction}

\begin{dfn}\label{metaunitarydatumdfn}
    A \emph{meta-unitary Shimura datum} is the data $\left(G,\mu,\rho_1,\rho_2,\rho_3\right)$
    where $(G,\mu)$ is a Shimura datum with $G = \Res_{L^+/\mathbb{Q}} G_0$, for some reductive group $G_0$ (see Condition 1 below) and for $i\in \{1,2,3\}$, $\rho_i$ is a skew-Hermitian representations of $G_0$ on Hermitian space $(V_i, \Psi_i)$ over $L^+$, satisfying the conditions below. Let $V_{123}$ denote the tensor product $V_1 \otimes_{L^+} V_2 \otimes_{L^+} V_3$.
    
    \begin{enumerate}        
        \item There exist involutive algebras $R_1, R_2, R_3$, and $R_{123}$, free and of finite rank over $\mathbb{Z}_{p^r}$, together with actions of $R_j$ on $V_j$ for $j \in \{1,2,3,123\}$, such that $G_0$ be the subgroup of $\prod_j \GL(V_j)$ consisting of those elements that commute with these actions;
        
        \item For $i\in \{1,2,3\}$, by Construction \ref{archemedianmutilde}, we obtain Hodge structures $(V^\triangleright_i, h_i^\triangleright)$. We assume that these Hodge structures satisfy Conditions \ref{hypothesiscompactness} and \ref{hypothesistensor};
        
        \item For $j\in \{1,2,3,123\}$ there exist a skew-Hermitian form $\Psi^\triangleright_j$ on $V^\triangleright_j$ making $(V^\triangleright_j, h_j^\triangleright, \Psi^\triangleright_j)$ is a skew-hermitian Hodge structure, and an isometry $\epsilon_j\colon V_j \otimes \mathbb{A}^\infty \otimes L \longrightarrow V^\triangleright_j \otimes \mathbb{A}^\infty \otimes L$;
        
        \item The datum is said to be \emph{unramified at $p$} if, moreover, there exists a self-dual $\mathcal{O}_{L,p}$-lattice $B_i \subset V_i$ such that $\textup{U}(B_i, \Psi_i)$ is hyperspecial inside $\textup{U}(V_i, \Psi_i)$;
    \end{enumerate}
    
    A meta-unitary Shimura datum is called \emph{regular} if the associated display datum is regular in the sense of Definition \ref{regular}.
\xqed{\blacktriangle}\end{dfn}

\begin{dfn}\label{dfnGtriangle}
We denote the group $\prod_j\Res_{L^+/\BQ} \textup{GU}(V_j,\Psi_j)$ by $G_\rho^\triangleright$.
\xqed{\blacktriangle}\end{dfn}

\begin{example}\label{example}
Let $\left(G,\mu\right)$ be a Shimura datum. Assume that $G$ splits over $\mathbb{Q}_{p^f}$, the derived subgroup $G^{\der}$ is simply connected, and the quotient adjoint group $G^{\ad}$ is simple and quasi-split over  $\mathbb{Q}_p$. Assume further that $G$ satisfies one of the two conditions:

    \begin{itemize}
        \item $G$ is of type $B_n$ or $C_n$ and $G^{\ad}\otimes \mathbb{R}$ has more than three times compact simple real factors than non-compact ones;
        \item $G$ is of type $E_7$ and $G^{\ad}\otimes \mathbb{R}$ has more than four times compact simple real factors than non-compact ones.
    \end{itemize}

Then, Appendix D of \cite{bultel2008pel} constructs an explicit meta-unitary Shimura datum for $\left(G,\mu\right)$. Some of these examples are of Abelian type, but there are also examples of type $E_7$ and of type $D_n$.

One can also construct regular meta-unitary Shimura data, using the fact that we can control which embeddings of $G$ are compact by the exact sequence, where $\pi_1(G)$ is the algebraic fundamental group of $G$ \cite[Proposition 2.6]{kottwitz1986elliptic}
    \begin{equation}
        H^1\left(L^+,G\right) \to \bigoplus_v H^1\left(L^+_{v},G\right) \to \pi_1\left(G\right)_{\Gal_{L^+}}.
    \end{equation}
\xqed{\blacksquare}\end{example}

Fix a compact open subgroup $K^p\subset G(\BA_{\mathbb{Q}}^{p,\infty})$. Consider the following PEL-type moduli problem associated with a meta-unitary Shimura datum as above:

\begin{construction}\label{hypothesisconstruction}
    Let $Y_1, Y_2$ and $Y_3$ be Abelian varieties over $\BC$, each equipped with an action of $\mathcal{O}_L$.
Write $H_i = H^1(Y_i(\BC), \BQ)$ for their rational Hodge structures endowed with the induced $\mathcal{O}_L$-actions, and let
\[
\boldsymbol{\mu}_i \colon \mathbb{G}_{m,\BC} \longrightarrow \BC^\times
\]
be the cocharacter associated with the Hodge structure $H_i$.
Under Condition \ref{hypothesistensor},  
the tensor product
\[
H \coloneqq H_1 \otimes_{\mathcal{O}_L} H_2 \otimes_{\mathcal{O}_L} H_3
\]
carries a natural Hodge structure of weight $1$, and hence corresponds to a complex Abelian variety with a canonical $\mathcal{O}_L$-action, called the \emph{tensor product Abelian variety} associated with $Y_1, Y_2$ and $Y_3$, and denoted by $Y_1\otimes Y_2\otimes Y_3$

Moreover, this construction extends functorially to families of Abelian varieties with $\mathcal{O}_L$-action over more general bases; see \cite{bultel2008pel} for details.
\xqed{\blacktriangledown}\end{construction}

\begin{dfn}
Let $R_1, R_2, R_3$, and $R_{123}$ be involutive $\mathcal{O}_L$-algebras that are free and of finite rank over $\mathcal{O}_L$.  

We define a moduli problem $\bar{\mathcal{S}}^{\textup{tan}}_{K^p}$ over $\mathbb{F}_q$ as follows.  
For every connected $\mathbb{F}_q$-scheme $X$,  
the set $\bar{\mathcal{S}}^{\textup{tan}}_{K^p}(X)$ consists of the data:
\begin{itemize}
    \item[(1)] A polarized Abelian scheme 
    \[
      (Y, \lambda, \eta, \iota)
    \]
    over $X$, equipped with
    \begin{itemize}
        \item a prime-to-$p$ level structure $\eta$ of type $K^p$,
        \item an action $\iota\colon \mathcal{O}_L^{\oplus 3} \to \End_X(Y)$,
        \item and a polarization $\lambda$
    \end{itemize}
    satisfying the usual PEL-type conditions.
    
    \item[(2)] For each $i \in \{1,2,3\}$, an action
    \[
      y_i \colon R_i \longrightarrow \End_S(Y_i), \qquad
      \text{where } Y_i \coloneqq e_i \cdot Y
    \]
    and $e_1, e_2$ and $e_3$ are the standard basis of $\mathcal{O}_L^{\oplus 3}$.
    These actions are required to satisfy the determinant condition relative to the given $\mathcal{O}_L$-structure.
    
    \item[(3)] An action of $R_{123}$ on 
    $Y_{123}\coloneqq Y_1\otimes Y_2\otimes Y_3$ (see Construction \ref{hypothesisconstruction}), compatible with its natural action on the prime-to-$p$ adelic Tate module.  
    More precisely, the diagram
    \[
    \begin{tikzcd}
    \mathbb{A}^{\infty,p}\! \otimes V_{123} \arrow{r}{\eta} &
    \displaystyle\bigotimes_{i=1}^3 H^1_{\textup{\'et}}(Y_{i,x_0}, \mathbb{A}^{\infty,p}) \arrow{r} &
    H^1_{\textup{\'et}}(Y_{123,x_0}, \mathbb{A}^{\infty,p})
    \end{tikzcd}
    \]
    commutes with the action of $R_{123}$.
\end{itemize}
\xqed{\blacktriangle}\end{dfn}

\begin{lem}\label{lem:closed-immersion}
    There is a natural forgetful morphism to the PEL-type Shimura variety
    \[
        f \colon \bar{\mathcal{S}}^{\textup{tan}}_{K^p} \longrightarrow \CS_{K^p}({G^\triangleright_\rho,\mu^\triangleright}),
    \]
    which is representable by a closed immersion. Consequently, $\bar{\mathcal{S}}^{\textup{tan}}_{K^p}$ is representable by a projective scheme.
\end{lem}

\begin{proof}
    The additional data parameterized by $\bar{\mathcal{S}}^{\textup{tan}}_{K^p}$, namely, the compatible actions of the  involutive algebras
    $R_1, R_2, R_3$ and $R_{123}$ on the corresponding factors of the Abelian variety, impose closed conditions
    on the moduli problem $\CS_{K^p}(G_\rho^\triangleright,\mu^\triangleright)$.
    Indeed, each required commutativity or determinant condition translates into the vanishing of certain morphisms between finite locally free sheaves.
    Therefore, the forgetful morphism $f$ is representable by a closed immersion.
\end{proof} 

\begin{rmk}
A priori, it is not even clear that this moduli space is non-empty. 
\xqed{\lozenge}\end{rmk}
In \cite{bultel2008pel}, B\"ultel constructed the functor $\Flex$ and proved that it is representable by a closed immersion. Now we can construct the candidate for the special fiber of $\Sh_{K}\left(G,\mu\right)$:

\begin{dfn}
    We denote the unique scheme that makes the following diagram Cartesian by $\bar{S}^{\textup{nr}}$.
        \[
        \begin{tikzcd}
        \bar{\CS}_{K}^{\textup{nr}} \arrow{r} \arrow{d} 
        & \CB^{\textup{Witt}}_{(G,\mu)} \arrow{d}{\Flex^{\textup{tan}}} \\
    \bar{\mathcal{S}}^{\textup{tan}}_{K^p} \arrow{r} 
        & \CB^{\textup{Witt},\textup{tan}}_{(G,\mu)}
        \end{tikzcd}
        \]
    We denote by $\bar{\CS}_{K}$ be the underlying reduced sub-scheme of $\bar{\CS}_{K}^{\textup{nr}}$.
\xqed{\blacktriangle}\end{dfn}

\begin{pro}\cite[Section 8.2.2]{bultel2008pel}
    The scheme $\bar{\CS}_K$ is smooth.
\end{pro}

\begin{proof}
    Note that $\CB^{\textup{Witt}}_{\left(G,\mu\right)}$ is a formally smooth stack \cite{bueltel2020g}, and the map $$\bar{\mathcal{S}}^{\textup{tan}}_{K^p}\to \CB^{\textup{Witt},\textup{tan}}$$ is formally \'etale by Serre-Tate, and so its base change $$\bar{\CS}_K^{\textup{nr}}\to \CB^{\textup{Witt}}_{\left(G,\mu\right)}$$ is also formally \'etale.
\end{proof}

\begin{pro}\label{PropLiftS-bar}
    There is a unique pair $\left(\CS_K, \CD\right)$, where $\CS_K$ is a smooth projective scheme over $W(\BF_q)$ that lifts $\bar{\CS}$, and $\CD$ is a $(G, \mu)$-display over $\CS_K$ that lifts the universal $\left(G, \mu\right)$-display on $\bar{\CS}$, given by the map $\bar{\CS}\to \CB^{\textup{Witt}}_{\left(G,\mu\right)}$.
\end{pro}

\begin{proof}
    We provide the sketch of the proof; for details, check \cite[Section 8]{bultel2008pel}. By the invariance of the \'etale site under nilpotent thickenings, one can lift $\bar{\CS}$ to a unique smooth projective variety $\CS_n$ over $\CB^{\textup{Witt}}_{\left(G,\mu\right),W_n(\BF_q)}$. By taking the limit, we obtain a smooth formal projective scheme $\CS$. One can construct an ample line bundle on each $\CS_n$, and therefore on their limit. By GAGA for formal schemes, we obtain a smooth lift $\CS$ of $\bar{S}$ with a map to $\CB^{\textup{Witt}}_{\left(G,\mu\right)}$.
\end{proof}

\begin{rmk}
 The $(G, \mu)$-display on $\CS_K$ induces an isocrystal with $G$-structure on $S$ equipped with an integrable (flat) connection. In particular, this gives rise to a $G$-torsor $\mathcal{P}$ with a flat connection over $\CS_K$, and the Hodge filtration associated with the $(G, \mu)$-display defines a morphism  $\mathcal{P} \longrightarrow \Gr_{G, {\CS}_K}$.
\xqed{\lozenge}\end{rmk}

In what follows, let $\CS_K$ be a smooth scheme over $\mathbb{Z}_q$ endowed with a $G$-torsor as above.  
To consider its complex fiber, we fix an algebraic isomorphism 
\[
    \iota\colon \overline{\mathbb{Q}}_p \xrightarrow{\sim} \mathbb{C},
\]
and define
\[
    S_{\mathbb{C}} \coloneqq S \otimes_{\mathbb{Z}_q,\,\iota} \mathbb{C}.
\]
We write $S_{\mathbb{C}}^{\mathrm{an}}$ for the associated complex analytic manifold and fix a base point $s_0 \in S_{\mathbb{C}}^{\mathrm{an}}$.

Let $\tilde{S}$ denote the universal cover of $S_{\mathbb{C}}^{\mathrm{an}}$.  
Define a group $\Delta$ whose elements are pairs $(g,h)$, where $g \in G(\mathbb{A}^\infty)$ 
and $h$ is the homotopy class of a path in $S_{\mathbb{C}}^{\mathrm{an}}$ from $s_0$ to $g \cdot s_0$.  
The group $\Delta$ acts on $\tilde{S}$, and one can identify
\[
    \Delta \backslash \bigl(\tilde{S} \times G(\mathbb{A}^{\infty,p})\bigr)
\]
with the union of certain connected components of $S_{\mathbb{C}}^{\mathrm{an}}$.  
For a point $y \in \tilde{S}$, let $\Gamma_y$ denote its stabilizer in $\Delta$.

The monodromy of the $G$-torsor with flat connection $\mathcal{P}$ gives rise to a representation
\[
    \tau\colon \Delta \longrightarrow G(\mathbb{C}).
\]

\begin{dfn}
    An \emph{elliptic point} of $\tilde{S}$ is a point $y \in \tilde{S}$ such that the stabilizer $\Gamma_y$ 
    contains a subgroup whose image under $\tau$ is an elliptic torus in $G(\mathbb{C})$.
\xqed{\blacktriangle}\end{dfn}

\medskip

Let us summarize the main results concerning $S$, as proved in \cite{bultel2008pel}.

\begin{thm}\label{mainproperties}
    There exists on $S$ a crystal with $G$-structure equipped with an integrable connection, 
    whose monodromy group is maximal. 
    Moreover, every connected component of $S$ contains an elliptic point.
\end{thm}

\begin{rmk}
    The proof of the existence of elliptic points is the only place in \cite{bultel2008pel}
    where the PEL-type assumption in the definition of meta-unitary Shimura varieties
    plays an essential role. 
    If this argument could be extended beyond the PEL setting—say, to Hodge-type data—then 
    one could similarly construct integral models in that broader context.
\xqed{\lozenge}\end{rmk}

\medskip

\begin{thm}\label{complex_shimura_variety}
Using the above properties and the characterization results of \cite{varshavsky2002characterization}, 
$S_{\mathbb{C},K^p}$ is isomorphic to the Shimura variety
    $\Sh(G,\mu)_{K^pK_p}$,
for a suitable compact open subgroup $K_p \subset G(\mathbb{Q}_p)$.
\end{thm}

\section{Shtukas and Local Shimura Varieties}\label{shtukaappaendix}

\subsection{Mixed-characteristic shtukas with several legs}

In this section, we recall the definition of local Shimura varieties. For more details, check \cite[Section 24]{scholze2020berkeley}.

Let $\CG$ be a reductive group over $\BZ_p$. Fix cocharacters of $\CG$, $\boldsymbol{\mu}=(\mu_i)_{i=1}^m$ where $\mu_i$ is defined over a finite extension $E_i$ of $\BQ_p$. Let $S$ be a perfectoid space and let
\[
x_1,\dots,x_m\colon S \longrightarrow \prod \mathrm{Spd}\,\BO_{E_i}
\]
be a collection of morphisms, called \emph{legs}. Write $\Gamma_{x_i}\subset \mathcal{Y}_S$ for the corresponding Cartier divisors.

\begin{dfn}\label{defshtukawithseverallegs}

A \emph{(mixed-characteristic) $\CG$-shtuka over $S$ with legs $x_1,\dots,x_m$ and bounded by $\boldsymbol{\mu}$} is given by:

\begin{itemize}
  \item a $\CG$-torsor $\mathcal E$ on $\mathcal{Y}_S$, together and
  \item a Frobenius-linear isomorphism
  \[
    \phi\colon\sigma^\ast(\mathcal E)\Big|_{\mathcal{Y}_S\setminus\bigcup_i\Gamma_{x_i}}
    \;\xrightarrow{\;\simeq\;}\;
    \mathcal E\Big|_{\mathcal{Y}_S\setminus\bigcup_i\Gamma_{x_i}},
  \]
  which is bounded by $\boldsymbol{\mu}$ along the divisors $\Gamma_{x_i}$. More precisely, the relative position of $\sigma^*\mathcal{E}$ and $\CE$ at $x_i$ is bounded by $\sum_{j|x_i=x_j} \mu_j$.
\end{itemize}
\xqed{\blacktriangle}\end{dfn}
We denote the moduli space of $\CG$-shtukas with $m$ legs over $\Spd \breve{E}_1\times \dots\times \Spd \breve{E}_m$ by $\Sht^m_{(G,\mu^\bullet)}$.

\subsection{Local Shimura varieties}

\begin{dfn}
    A \emph{local Shimura datum} is a triple $(G,\{\mu\},[b])$ consisting of:
    \begin{itemize}
      \item a connected reductive group $G$ over $\BQ_p$,
      \item a conjugacy class of cocharacters $\{\mu\}$ of $G$ defined over some extension $E/\BQ_p$,
      \item A $\sigma$-conjugacy class $[b]\subset G(\breve{\BQ}_p)$.
    \end{itemize}
    By abuse of notation we denote such a triple by $(G,\mu,b)$.
\xqed{\blacktriangle}\end{dfn}

Scholze and Weinstein associate with $(G,\mu,b)$ a tower of diamonds called the \emph{local Shimura variety} associated with $(G,\mu,b)$ 

\begin{dfn}[Local Shimura variety]
For a compact open subgroup $K\subset G(\BQ_p)$, the \emph{local Shimura variety}
\[
\mathcal{M}_{(G,\mu,b),K}
\]
is the diamond over $\mathrm{Spd}\,\breve E$ representing the following functor:
for a perfectoid $\breve E$-space $S$, the groupoid
$\mathcal{M}_{(G,\mu,b),K}(S)$ consists of
\begin{itemize}
  \item a $G$-torsor $\mathcal E$ on the relative Fargues--Fontaine curve $X_S$,
  \item together with an isomorphism
  \[
    \alpha\colon\; \mathcal E\big|_{X_S\setminus\{x\}}
    \;\xrightarrow{\ \simeq\ }\;
    \mathcal E_b\big|_{X_S\setminus\{x\}},
  \]
  for a leg $x\colon S \to \mathrm{Spd}\,\BQ_p$,

  \item a $K$-lattice inside the $G(\mathbb{Q}_p)$-torsor associated to $\mathcal{E}$.
\end{itemize}
such that the modification $(\mathcal E,\mathcal E_b,\alpha)$ is bounded by $\mu$, i.e.\ its image in the affine Grassmannian $\mathrm{Gr}_G$ lies in the Schubert variety $\mathrm{Gr}_G^{\leq \mu}$.
\xqed{\blacktriangle}\end{dfn}

This tower is representable by a tower of rigid spaces. The tower carries a natural action of $G(\BQ_p)$ (via Hecke correspondences) and of the group $J_b(\BQ_p)$ of self quasi-isogenies of the isocrystal attached to $b$.

\subsection{Integral local Shimura varieties}
Let $(\CG,\{\mu\},[b])$ be an integral local Shimura datum: $\CG$ a connected reductive group over $\mathbb{Z}_p$, $[b]\subset \CG(\breve\BQ_p)$ a $\sigma$-conjugacy class, and $\{\mu\}$ a conjugacy class of cocharacters of $\CG$ defined over a finite extension $E$ of $\BQ_p$.   

For a perfectoid space $S$ over $\Spd(\CO_{\breve E})$, denote by $\mathcal E_b$ the $\CG$-bundle on $X_S$ corresponding to a representative $b\in \CG(\breve\BQ_p)$.

\begin{dfn}\label{dfnintegrallocal}
The integral local Shimura variety associated with $(\CG,\mu,b)$ is defined as the functor $\mathcal{M}^{\textup{int}}_{(\CG,\mu,b)}$ sending $S\overset{x}{\to}\Spd \CO_{\breve{E}}$ to the set of isomorphism classes of triples $(\mathcal E,\;F,\;\eta) $ where

\begin{itemize}
    \item $\CE$ is a $\CG$-torsor on $\mathcal{Y}_S$,
    \item $F\colon \sigma^*\mathcal E\xrightarrow{\simeq}\mathcal E$ is a modification at the leg $x\colon S\to\Spd\CO_{\breve{E}}$ bounded by $\mu$, and
    \item $\eta\colon\mathcal{E}_{[r,\infty)}\xrightarrow{\simeq}\mathcal{E}_{b,[r,\infty)}$ is an isomorphism.
\end{itemize}
\xqed{\blacktriangle}\end{dfn}

\end{appendices}
\small
\bibliography{myref}{}
\bibliographystyle{alphaurl}

\end{document}